\documentclass[10pt]{amsart}
\usepackage{setspace, pb-diagram, lamsarrow,pb-lams}
\usepackage{amsmath, amscd, amsthm, amssymb, amsfonts}
\usepackage{verbatim, array, longtable}
\title[Gromov--Witten invariants and symmetric products]{Higher genus Gromov--Witten invariants as genus zero invariants of symmetric
products}
\author{Kevin Costello}
\address{
Department of Mathematics \\  Imperial College \\ London  SW7 2AZ\\
United Kingdom}
 \email
{k.costello@imperial.ac.uk}
\date{}

\newcommand{\iso}{\cong}
\newcommand{\diag}{\bigtriangleup} 
\newcommand{\C}{\mathbb C}

\newcommand{\norm}[1]{\left\| #1 \right\|}
\newcommand{\Oo}{\mathcal O}
\newcommand{\fock}{\mathcal F}
\newcommand{\Z}{\mathbb Z}
\newcommand{\defeq}{\overset{\text{def}}{=}}
\newcommand{\into}{\hookrightarrow}

\newcommand{\op}{\operatorname}

\newcommand{\mbb}{\mathbb}
\newcommand{\mf}{\mathfrak}
\newcommand{\mc}{\mathcal}
\newcommand{\from}{\leftarrow}
\newcommand{\ip}[1]{\left\langle #1 \right\rangle}

\newcommand{\vertex}{\Upsilon^t}
\newcommand{\rvertex}{\Upsilon^u}
\newcommand{\cvertex}{\Upsilon^c}
\renewcommand{\mod}{\mathfrak M}
\newcommand{\cover}{\twoheadrightarrow}
\newcommand{\curve}{\mathfrak C}
\newcommand{\graph}{\Gamma^t}
\newcommand{\rgraph}{\Gamma^u}
\newcommand{\cov}{\Gamma^c}
\newcommand{\cont}{\rightarrow}

\newcommand{\cmod}{\overline{ \mathcal M}}

\newcommand{\Q}{\mbb Q}

\DeclareMathOperator{\Aut}{Aut}

\DeclareMathOperator{\Sym}{Sym} \DeclareMathOperator{\Hom}{Hom}
\DeclareMathOperator{\Spec}{Spec}

\newtheorem{theorem}{Theorem}[subsection]
\newtheorem{proposition}[theorem]{Proposition}

\newtheorem{lemma}[theorem]{Lemma}

\newtheorem{corollary}[theorem]{Corollary}

\newtheorem{example}[theorem]{Example}

\numberwithin{equation}{subsection}

\begin{document}

\begin{abstract}
I prove a formula expressing the descendent genus $g$ Gromov-Witten invariants
of a projective variety $X$ in terms of genus $0$ invariants of its symmetric
product stack $S^{g+1}(X)$. When $X$ is a point, the latter are structure
constants of the symmetric group, and we obtain a new way of calculating the
Gromov-Witten invariants of a point.
\end{abstract}

\maketitle
\section{Introduction}
Let $X$ be a smooth projective variety.  The genus $0$ Gromov-Witten invariants
of $X$ satisfy relations which imply that they can be completely encoded in the
structure of a Frobenius manifold on the cohomology $H^\ast(X,\C)$.     In this
paper I prove a formula which expresses the descendent genus $g$ Gromov-Witten
invariants of a smooth projective variety $X$ in terms of the descendent genus
$0$ invariants of the symmetric product stack $S^{g+1}X$. The latter are
encoded in a Frobenius manifold structure on the orbifold cohomology group
$H^\ast_{orb}(S^{g+1}(X),\C)$. This implies that the Gromov-Witten invariants
of $X$ at all genera are described by a sequence of Frobenius manifold
structures on the homogeneous components of the Fock space
\begin{equation*}
  \fock = \Sym^\ast \left( H^{\ast}(X, \C) \otimes_\C t \C[t] \right) =
  \oplus_{d \ge 0} H^\ast_{orb} ( S^d (X), \C)
\end{equation*}
Standard properties of genus $0$ invariants, such as associativity, when
applied to the symmetric product stacks $S^d X$, yield implicit relations among
higher-genus Gromov-Witten invariants of $X$.

When $X = \ast$ is a point, the symmetric product is the classifying stack
$BS_d$ of the symmetric group. The Frobenius manifold associated to the genus
$0$ invariants of $BS_d$ is in fact a Frobenius \emph{algebra}, which is the
centre of the group algebra of the symmetric group, $\C[S_d]^{S_d}$. Our result
therefore gives a new way of expressing the integrals of tautological classes
on $\cmod_{g,n}$ in terms of structure constants of $\C[S_d]$.

More generally, the associativity constraints, together with some other simple
properties, are sufficient to determine the small quantum cohomology of the
symmetric product stack $S^d X$ in terms of the small quantum cohomology of
$X$. The construction of Lehn-Sorger \cite{lehn-sorger}, (see also
Fantechi-G\"ottsche \cite{fantechi gottsche}), which calculates the orbifold
cohomology of $S^d X$ in terms of the ordinary cohomology of $X$, applies
verbatim to calculate the small quantum cohomology of $S^d X$ in terms of that
of $X$. In general, the large quantum cohomology of $S^d X$ is not determined
by that of $X$.

Let me sketch the geometric relation between Gromov-Witten invariants of $X$
and $S^d X$. Stacks of stable maps to the symmetric product stack $S^d X$ are
identified with stacks of certain correspondences $\mc C \from \mc C' \to X$,
where $\mc C$ and $\mc C'$ are twisted balanced nodal curves, and $\mc C' \to
\mc C$ is \'{e}tale of degree $d$. We introduce stacks $\cmod_\eta(X)$,
parameterizing such correspondences with certain markings on $\mc C'$ and $\mc
C$, where $g(\mc C) = 0$. \footnote{ We identify the genus of a twisted curve
with that of its coarse moduli space.} Here $\eta$ is some label remembering
the genera of $\mc C'$ and $\mc C$, the stack structure at the marked points,
the homology class of the map $\mc C' \to X$, and so forth. There is a finite
group $G$ acting without fixed points on $\cmod_\eta(X)$, by reordering marked
points of $\mc C'$, such that $\cmod_\eta(X)/G$ is a stack of stable maps from
genus $0$ curves to $S^d X$. This implies that integrals on $\cmod_\eta(X)$ are
Gromov-Witten invariants of $S^d X$.

There is a map $p: \cmod_\eta(X) \to \cmod_{g,r,\beta}(X)$, for some $g,r$ and
$\beta \in H_2(X)$, defined by taking the coarse moduli space $C'$ of $\mc C'$,
with its natural map  $C' \to X$, and forgetting some marked points.  We show
that $p$ is finite of degree $k \in \Q^\times$, in the virtual sense.  By this
we mean
\begin{equation}
  p_\ast [\cmod_{\eta}(X)]_{virt} = k [\cmod_{g,r}(X) ]_{virt}
  \label{intro_push_forward}
\end{equation}
We then express the pull back $p^\ast \psi_i$ of the tautological $\psi$
classes on $\cmod_{g,r}(X)$, in terms of $\psi$ classes and boundary divisors
of $\cmod_{v' \to v}( X)$.  The boundary cycles of $\cmod_{\eta}(X)$ are again
products of similar stacks of \'{e}tale correspondences.  Further,  there is a
commutative diagram of evaluation maps
\begin{equation*}
  \begin{diagram}
    \node{\cmod_{\eta} (X)} \arrow{e} \arrow{s} \node{X^r} \\
    \node{\cmod_{g,r}(X)} \arrow{ne}
  \end{diagram}
\end{equation*}
This allows us to translate integrals on $\cmod_{g,r}(X)$ of $\psi$ classes,
and cohomology classes pulled back from $X^r$, into sums of products of similar
integrals on $\cmod_{0,n}(S^m X)$ for varying $m$ and $n$.

The most technically difficult part of this procedure is proving the
push-forward formula (\ref{intro_push_forward}).  We do this by working in a
``universal'' setting, where all the moduli stacks are smooth; and deduce it
for arbitrary $X$ by base change, in the virtual sense, by the stack
$\cmod_{g,\beta}(X)$ of curves in $X$ with no markings.  We need to introduce
moduli stacks of curves with markings in a semigroup. Let $A$ be a semigroup
with indecomposable zero; for each $a \in A$ we define a moduli stack
$\mod_{g,n,a}$ of all (possibly unstable) connected nodal curves of genus $g$,
with $n$ marked points, and certain $A$-valued marking on the irreducible
components.  These curves must satisfy some stability conditions; for example
when $a = 0$, but not otherwise, $\mod_{g,n,0} = \cmod_{g,n}$ is the usual
Deligne-Mumford moduli stack of stable curves. In general, $\mod_{g,n,a}$ is
smooth, proper, locally of finite type, but non-separated. The advantage of
these moduli stacks over the more familiar stacks $\mod_{g,n}$ of all nodal
curves, is that there are (proper, separated) contraction maps $\mod_{g,n,a}
\to \mod_{g,n-1,a}$, which identify $\mod_{g,n,a}$ with the universal curve
over $\mod_{g,n-1,a}$. This is not the case for $\mod_{g,n}$.

Let $C(X)$ be the Mori cone of positive 1-cycles on $X$ modulo numerical
equivalence. For each $\beta \in C(X)$ we have the associated smooth moduli
stacks $\mod_{g,n,\beta}$, and the stacks of stable maps
$\cmod_{g,n,\beta}(X)$.  We have
\begin{equation*}
  \cmod_{g,n,\beta}(X) = \mod_{g,n,\beta} \times_{\mod_{g,\beta}}
  \cmod_{g,\beta}(X)
\end{equation*}
More generally, for any connected modular graph $\gamma$ with labellings in
$C(X)$, so that $\gamma$ defines a stratum of $\cmod_{g,n,\beta}(X)$, we see
that
\begin{equation*}
  \cmod_\gamma (X) = \mod_{\gamma} \times_{\mod_{g,\beta} }
  \cmod_{g,\beta}(X)
\end{equation*}
Further, these fibre products are compatible with virtual fundamental classes.
That is, the system of stacks of stable maps to $X$, together with their
natural morphisms and virtual classes, is pulled back, via the map
$\cmod_{g,\beta}(X) \to \mod_{g,\beta}$, from the stacks $\mod_{g,n,a}$ with
their natural morphisms.

We can extend this observation to stacks of \'{e}tale correspondences to $X$:
\begin{equation*}
 \cmod_{\eta}(X) = \mod_{\eta} \times_{\mod_{g,\beta}} \cmod_{g,\beta}(X)
\end{equation*}
where $\mod_\eta$ is some stack of \'{e}tale maps of curves $\mc C' \to \mc C$.
This fibre product is also compatible with virtual fundamental classes.

Since all of these fibre products are in the virtual sense, they behave quite
like flat base changes for the purposes of intersection theory. We show that to
prove the map $\cmod_{\eta}(X) \to \cmod_{g,r,\beta}(X)$ is finite in the
virtual sense as in formula (\ref{intro_push_forward}), it is sufficient to
show that $\mod_{\eta} \to \mod_{g,r,\beta}$ is actually finite.  With the
correct choices of $\eta$, this is not difficult.

\subsection{Relation to previous work}
\subsubsection*{Intersection numbers on moduli stacks of curves}
The Gromov-Witten theory of a point has been known since Kontsevich's proof
\cite{kontsevich} of Witten's conjecture \cite{witten}.  There are two parts to
Kontsevich's proof. Firstly, he reduces the geometric problem to a
combinatorial problem,  using a topological cell decomposition of the moduli
stack of curves to derive formulae for integrals of tautological classes. Then
he derives a matrix integral formula for these expressions, and uses this to
prove Witten's conjecture.

The results of this paper, applied to a point, give a new way to do the first
part of this procedure; that is we find a combinatorial expression for
integrals of tautological classes on the moduli stack. The techniques are
purely algebro-geometric, and thus have a very different flavour from
Kontsevich's topological model.

More recently, another proof of the Kontsevich-Witten theorem has appeared. A
combinatorial expression for intersection numbers on the moduli stack of curves
in terms of Hurwitz numbers was announced by Ekedahl, Lando, Shapiro and
Vainshtein in \cite{elsv1} and proved in \cite{elsv2}.  Another proof of this
formula was obtained by Graber and Vakil \cite{graber vakil}, building on a
special case proved by Fantechi and Pandharipande \cite{fantechi
pandharipande}. This result was used by Okounkov and Pandharipande \cite{ok
pan} to give another proof of the Kontsevich-Witten theorem.

The geometric part of this proof relies on spaces of ramified covers of $\mbb
P^1$ to relate intersection numbers on $\cmod_{g,n}$ to Hurwitz numbers. Spaces
of covers of genus $0$ curves also play a central role in this work. However,
the compactifications we use are different, as are the methods for obtaining
 formulae for integrals on $\cmod_{g,n}$.  For example, in
\cite{graber vakil}, Graber and Vakil calculate certain Gromov-Witten
invariants of $\mbb P^1$ in two different ways: firstly, using virtual
localization, and secondly, by using a ``branching map''  to configuration
spaces of points on $\mbb P^1$. Equating these yields the desired formula. On
the other hand,the techniques used here can be viewed, in the case of a point,
as firstly constructing a correspondence $\cmod_{0,n} \from \cmod_\eta \to
\cmod_{g,m}$, which is finite over both $\cmod_{0,n}$ and $\cmod_{g,m}$, and
then calculating the pullbacks of tautological classes from $\cmod_{g,m}$. The
expressions we end up with are different from those obtained by the authors
cited above.

The results presented here work for arbitrary target space, and not just a
point; it does not seem to be clear how to generalize the results of
\cite{elsv1,elsv2,fantechi pandharipande,graber vakil, kontsevich} to arbitrary
target $X$.

\subsubsection*{Orbifold Gromov-Witten theory}
Gromov-Witten invariants for orbifolds were first defined by Chen and Ruan
\cite{chen ruan, chen ruan2, chen ruan3}. In these papers they
introduced the orbifold cohomology groups $H^\ast_{orb}$, as the
ordinary cohomology of the space of twisted sectors of an orbifold. In
orbifold Gromov-Witten theory the orbifold cohomology group
$H^\ast_{orb}$ plays the same role as the ordinary cohomology group
plays in standard Gromov-Witten theory.  In particular, orbifold
quantum cohomology ($g = 0$ orbifold Gromov-Witten theory) gives
$H^\ast_{orb}$ the structure of a Frobenius manifold.  

Chen-Ruan's theory uses differential and symplectic geometry.
In algebraic geometry, Abramovich and Vistoli \cite{av} 
defined stable maps to Deligne-Mumford stacks,  and proved these form reasonable
stacks. In \cite{agv}, Abramovich, Graber and Vistoli use these results to give
an algebraic definition of Gromov-Witten invariants for DM stacks.

The Gromov-Witten theory of the classifying stack $BG$ of a finite group $G$
was studied by Jarvis and Kimura \cite{jarvis kimura}.    In a recent preprint,
Jarvis, Kaufmann and Kimura \cite{jarvis kaufmann kimura} study the algebraic
structure defined by $G$-equivariant quantum cohomology for a finite group $G$.
The reader should refer to these works for more details on the structure of the
genus $0$ Gromov-Witten invariants of $BS_n$ and of $S^n X$. Note, however,
that the notation for tautological classes, etc., used here, differs from their
notation by constants.

\subsection{Plan of the paper}

We define some of the basic moduli stacks we need in section
\ref{sec_moduli_spaces}. These are certain stacks of nodal curves with markings
in a semigroup; we show they are smooth Artin algebraic stacks and describe
certain maps between them, as well as tautological classes. Section
\ref{sec_graphs} is devoted to setting up various categories of labelled
graphs, together with functors which associate to a graph a certain moduli
stack of curves. In section \ref{sec_pull_back} we calculate how tautological
classes and cycles pull back under morphisms of moduli stacks, coming from
morphisms of graphs. These pull backs are expressed as sums over graphs,
weighted by tautological classes.

Section \ref{sec virtual} contains the main technical theorem, which says
roughly that a map of finite degree remains of finite degree in the virtual
sense, after a virtual base change.    We use Behrend-Fantechi's virtual
fundamental class technology, and this result follows from an analysis of their
``relative intrinsic normal cone stacks''.  In section \ref{sec_push_forward},
we construct, for each $g,r,\beta$, a label $\eta$ for a moduli stack of
\'{e}tale covers $\mod_\eta$, with a finite map $\mod_\eta \to
\mod_{g,r,\beta}$.

In section \ref{sec_base_change}, we base change by $\cmod_{g,\beta}(X)$, to
get stacks of stable maps and \'{e}tale correspondences to $X$. In section
\ref{sec main theorem}, we put these results together to give a formula for all
descendent Gromov-Witten invariants of $X$ in terms of genus $0$ invariants of
symmetric products $S^d X$. Finally, in section \ref{section calculations}, I
illustrate the general result by calculating some low-genus Gromov-Witten
invariants of a point.

\subsection{Future work}
The results presented here provide implicit constraints on the Gromov-Witten
invariants of an arbitrary variety, coming from associativity properties of
quantum cohomology of symmetric products.  It would be interesting to see what
relation these constraints have with the conjectural Virasoro constraints,
first proposed by Eguchi, Hori and Xiong \cite{eguchi hori xiong}. A first step
in this direction would be to use the results of this paper to give a new proof
of Witten's conjecture. I imagine this is far from easy; in the two proofs of
the Kontsevich-Witten theorem of which I am aware, even once the geometric work
has been done,  significant insight is required to prove the theorem.

In another direction, I think that one can prove a reconstruction theorem,
analogous to the first reconstruction theorem of Kontsevich and Manin
\cite{kontsevich manin}, which would imply that for certain Fano manifolds $X$,
the quantum cohomology of $S^d X$ is determined by the quantum cohomology of
$X$.  This would imply that all Gromov-Witten invariants of $X$ are determined
by the genus $0$ invariants. It's not clear in what generality one can make
such a statement: we need $K_X \ll 0$, which implies many genus $0$ invariants
of $S^d X$ vanish for dimension reasons.

Note that such a statement has a close relationship with certain corollaries of
the Virasoro conjecture. Dubrovin and Zhang \cite{dubrovin zhang 1, dubrovin
zhang 2} have shown that the Virasoro conjecture implies that when $X$ has
semi-simple quantum cohomology, all higher genus Gromov-Witten invariants of
$X$ are determined by the genus $0$ invariants. Conjecturally, many Fano
manifolds have semisimple quantum cohomology.

I hope to return to these points in a future paper.
\subsection{Acknowledgements}
I am very grateful to my Ph.D. supervisor Ian Grojnowski  for his support, both
mathematical and moral, over the last four years.  I'd like to thank Alessio
Corti, Ezra Getzler, Tim Perutz, Richard Thomas,  Burt Totaro, Hsian-Hua Tseng,
and especially Constantin Teleman and Ravi Vakil for their interest in this
work, and for very helpful conversations and correspondence. This paper will
form part of my Ph.D. thesis at Cambridge University. I have been financially
supported by the EPSRC, the Cambridge European Trust, and the Cecil King
Memorial Foundation.

\subsection{Notation}
We work always over a field $k$, algebraically closed and of characteristic
zero.  Stacks are in the sense of Laumon and Moret-Bailly \cite{lmb}. In
particular, a stack is not required to have an atlas, an algebraic stack must
admit a smooth atlas, and a Deligne-Mumford stack must admit an \'etale atlas.
I will sometimes use the phrase Artin stack as a synonym for algebraic stack.

Later we will define various categories of graphs. Here is a summary of some
notation needed for these:

\begin{tabular}{p{1.5in} p{3in} <{\vspace{5pt}} }
$\rvertex$ \hspace{2in} & Labels $(g,I,a)$ for smooth connected curves,
of genus $g$ with marked point  set $I$ and class $a \in A$ in the semigroup.\\
$\vertex$ & Labels $(g,I,m,a)$ for smooth connected twisted curves, of genus
$g$, with marked point set $I$, stack structure at the marked points given by
$m : I \to \Z_{> 0}$,  and class $a \in A$ in the semigroup .\\
$\op{exp}(\rvertex)$ \par $\op{exp}(\vertex)$ & Labels for disconnected smooth
marked curves, and disconnected smooth marked twisted curves, respectively.  \\
  $\cvertex$ & Labels for \'{e}tale covers $\mc C' \to \mc C$ of smooth   twisted
  curves, with $\mc C$ connected.\\
  $s : \cvertex \to \op{exp} (\vertex)$ & Source map, which associates to a
  label for $\mc C' \to \mc C$ the label for $\mc C'$. \\
  $t : \cvertex \to \vertex$ & Target map, which associates to a label for $\mc
  C' \to \mc C$ the label for $\mc C$. \\
  $\rgraph$ & Graphs built from vertices $\rvertex$, which label nodal
  connected curves.\\
  $\graph$ & Graphs built from vertices $\vertex$, which label  twisted
  nodal curves.\\
  $\cov$ & A certain type of map of graphs in $\graph$, which labels \'{e}tale
  covers $\mc C' \to \mc C$ of twisted nodal curves. \\
  $r: \vertex \to \rvertex$  \par $r: \graph \to \rgraph$ & Associates to a label for a twisted curve $\mc C$
  the label for its coarse moduli space $C$. \\
  $s,t: \cov \to \graph$ & Source and target maps, which associate to a label
  for an \'{e}tale cover $\mc C' \to \mc C$ the labels for $\mc C'$, $\mc C$
  respectively.
\end{tabular}

\section{Moduli stacks}
\label{sec_moduli_spaces}
Let $g \in \Z_{\ge 0}$ and let $I$ be a finite set.  Let $\mod_{g,I}$ be the
stack of all nodal curves of genus $g$ with $I$ marked smooth points.
$\mod_{g,I}$ is a smooth  algebraic stack; it is non-separated, and locally but
not globally of finite type.

Let $g \in \Z_{\ge 0}$, let $I$ be a finite set and let $m: I \to \Z_{> 0}$ be
a function.  Let $\mod_{g,I,m}$ be the moduli stack of all twisted (balanced)
curves of genus $g$, with  marked points labelled by $I$ and the degree of
twisting at the marked points given by $m$. Explicitly, $\mod_{g,I,m}$ is the
category of commutative diagrams
\begin{equation*}
  \begin{diagram}
     \node[2]{\mc C} \arrow{e} \arrow{se} \node{C} \arrow{s}\\
   \node{U \times I} \arrow{e}    \node{U \times \coprod_{i \in I} B\mu_{m(i)}} \arrow{n,J} \arrow{e} \node{U}
  \end{diagram}
\end{equation*}
where:
\begin{itemize}
\item
$U$ is a scheme of finite type.
\item
    $\mc C$ is a  proper separated flat DM stack over $U$, \'{e}tale locally a nodal curve over $U$.
\item
    The map $\mc C \to C$ exhibits $C$ as the coarse moduli space of $\mc C$, and $C$
    is connected of genus $g$.
\item
    $U \times \coprod B \mu_{m(i)} \into \mc C$ is an embedding of a disjoint
    union of trivial $\mu_{m(i)}$-gerbes into $\mc C$, and $U \times I \to U
    \times \coprod B\mu_{m(i)}$ are sections of these gerbes.
\item
    $\mc C \to C$ is an isomorphism away from the nodes and marked points of
    $C$.
\item
    \'Etale locally near a node of $\mc C$, $\mc C \to U$ looks like
\begin{equation*}
  (\Spec A[u,v]/ (uv - t))/ \mu_r \to \Spec A
\end{equation*}
    where $ t \in A$, and the group of $r$-th roots of unity $\mu_r$ acts on
    $A[u,v]/(uv-t)$ by $u \to lu$, $v \to l^{-1}v$,  where $l \in \mu_r$.
\end{itemize}
This definition is due to Abramovich and Vistoli; for more details see
\cite{av}. Note that we use trivialized gerbes, where they use possibly
non-trivial gerbes. Our stack is simply the fiber product of all the universal
gerbes lying over their version.
\begin{proposition}
$\mod_{g,I,m}$ is a smooth stack.
\end{proposition}
By smooth I mean in the sense of the formal criterion for smoothness over the
base $\op{Spec} k$.   I expect that $\mod_{g,I,m}$ is algebraic, although I
don't know a reference for this. Presumably, one could prove this using the
techniques of Abramovich and Vistoli \cite{av}. However, we don't really need
any properties of $\mod_{g,I,m}$; for us it is essentially a placeholder.

There is a map $\mod_{g,I,m} \to \mod_{g,I}$ which associates to a twisted
curve its coarse moduli space.

We need variants of these definitions, which depend on a semigroup.  Let $A$ be
a commutative semigroup, with unit $0 \in A$, such that
\begin{itemize}
  \item $A$ has indecomposable zero: $a + a' = 0$ implies $a = a' = 0$.
  \item $A$ has finite decomposition: for every $a \in A$, the set $\{ (a_1,
  a_2) \in A \times A \mid a_1 + a_2 = a \}$ is finite.
\end{itemize}
For example, $A = 0$, or $A$ is the Mori cone $C(X)$ of curves in a projective
variety $X$ up to numerical equivalence, or $A = \{0,1\}$ where $1 + 1 = 1$.

Fix any such $A$. Let $(g,I,a)$ be a triple where $g \in \Z_{\ge 0}$, $I$ is a
finite set, and $a \in A$. We say $(g,I,a)$  is stable, if either $a \neq 0$ or
$a = 0$ and $2g-2 + \#I
> 0$.  For any such triple $(g,I,a)$ we define the stack $\mod_{g,I,a}$ over
$\mod_{g,I}$.  Roughly, $\mod_{g,I,a}$ parameterizes curves $C$ with $I$ marked
smooth points, together with a labelling of each irreducible component
 of $C$ by an element of $A$.  The sum over irreducible components of the associated elements of $A$ must be
$a$, and a certain stability condition must be satisfied.  The simplest formal
definition is inductive.
\begin{enumerate}
  \item If $(g,I,a)$ is unstable, then $\mod_{g,I,a}$ is empty.
  \item Suppose $(g,I,a)$ is stable. Then an object of $\mod_{g,I,a}$ is
  \begin{itemize}
   \item An object of $\mod_{g,I}$,
   that is, a flat family $C \to U$, of nodal
  curves over a scheme $U$, together with $I$ smooth marked points $U \times I
  \to C$.
  \item
  Let $C_{\op{gen}} \to U$ be the complement of the nodes and marked points of $C$.
  The additional data we require is a constructible function $f: C_{\op{gen}} \to A$.
  $f$ must be locally constant on the geometric fibres of $C_{\op{gen}} \to U$.
  \end{itemize}
  \item If $U_0 \subset U$ is the open subscheme parameterizing non-singular curves $C_0 \to U_0$, then
   $f: {C_0}_{\op{gen}} \to A$ must be constant with value $a$.
  \item
   We require that $f$ satisfies a gluing condition along the boundary of
   $\mod_{g,I}$.  Precisely, suppose we have a decomposition $g = g' + g''$ and
   $I = I' \coprod I''$,  a map $V \to U$, and a
   factorization of the map $V \to \mod_{g,I}$ into
\begin{equation*}
  V \to \mod_{g',I'\coprod\{s'\}} \times
   \mod_{g'',I''\coprod \{s''\}} \to \mod_{g,I}
\end{equation*}
    where the second map is obtained by gluing the marked points $s', s''$.
    Let $C'_V \to V$ and $C''_V \to V$ be the associated families of curves.
    We require that the pulled-back constructible functions $f' : {C'_V}_{\op{gen}} \to A$
    and $f'': {C''_V}_{\op{gen}} \to A$ define a morphism
\begin{equation*}
  V \to \coprod_{a = a' + a''} \mod_{g',I'\coprod\{s'\},a'} \times
   \mod_{g'',I''\coprod \{s''\},a''}
\end{equation*}
\item
In a similar way, suppose we have a map $V \to U$, and a factorization of the
map $V \to \mod_{g,I}$ into
\begin{equation*}
  V \to \mod_{g-1,I \coprod \{s,s'\}} \to \mod_{g,I}
\end{equation*}
Then, the family of genus $g-1$ curves $C_V \to V$, together with the
pulled-back constructible function $f: {C_V}_{\op{gen}} \to A$, must define a
map
\begin{equation*}
  V \to \mod_{g-1,I \coprod\{s,s'\},a}
\end{equation*}
\end{enumerate}
\begin{proposition}
The map $\mod_{g,I,a} \to \mod_{g,I}$ is  \'{e}tale, and relatively a scheme of
finite type.   Therefore $\mod_{g,I,a}$ is a smooth algebraic stack.
\end{proposition}

Define $\mod_{g,I,m,a} = \mod_{g,I,m} \times_{\mod_{g,I}} \mod_{g,I,a}$.  The
stack $\mod_{g,I,m,a}$ is smooth.

\subsection{Contraction maps}
The main advantage of $\mod_{g,I,a}$ over $\mod_{g,I}$ is the existence of
contraction maps $\pi_i : \mod_{g,I,a} \to \mod_{g,I \setminus i,a}$ for each
$i \in I$. Given a curve $C \in \mod_{g,I,a}$ with  marked points $P_j$, $j \in
I$, $\pi_i(C)$ is obtained from $C$ by removing $P_i$, and contracting the
irreducible component of $C$ containing $P_i$ to a point if it is unstable.
Unstable components are components of genus $0$ with marking $0 \in A$ and
containing $< 3$ nodes or marked points, and components of genus $1$ with no
nodes or marked points. To construct the map $\pi_i$, we need
\begin{proposition}
There is an isomorphism $\curve_{g,I \setminus i, a} \iso \mod_{g,I,a}$ where
$\curve_{g,I,\setminus i,a}$ is the universal curve over $\mod_{g,I\setminus
i,a}$.
\end{proposition}
\begin{proof}
Just as in \cite{knudsen}, Definition 2.3, there is a map $\curve_{g,n-1} \to
\mod_{g,n}$.  This lifts to a map $\curve_{g,n-1,a} \to \mod_{g,n,a}$, by
labelling any irreducible component which is contracted in the map
$\curve_{g,n-1,a} \to \mod_{g,n-1,a}$ by $0 \in A$.  As in \cite{behrend gw},
lemma 7, the formal criterion for etaleness shows that $\curve_{g,n-1,a} \to
\mod_{g,n,a}$ is \'{e}tale.  To show it is an isomorphism, it is enough to show
this on the level of $k$-points, which is easy.
\end{proof}

\subsection{Maps to symmetric products}
\label{sec map sym prod}
Let $X$ be a scheme. Let $S^d X = [ X^d / S_d]$ be the symmetric product stack
of $X$.
\begin{lemma}
The stack whose groupoid of $U$ points, for $U$ a scheme, has objects diagrams
\begin{equation*}
  U \from U' \to X
\end{equation*}
where $U' \to U$ is proper, separated, surjective, and \'{e}tale of degree $d$;
and has morphisms, isomorphisms $U' \to U'$ such that the diagram
\begin{equation*}
  \begin{diagram}
    \node[2]{U'} \arrow{sw} \arrow[2]{s} \arrow{se} \\
    \node{U} \node[2]{X} \\
    \node[2]{U'} \arrow{nw} \arrow{ne}
  \end{diagram}
\end{equation*}
commutes, is equivalent to the stack $S^d X$.
\end{lemma}
\begin{proof}
By definition, to give a map $U \to S^d X$ is to give a right, \'{e}tale
locally trivial,  principal $S_d$-bundle $P \to U$, together with an
$S_d$-equivariant map $P \to X^d$.  Given such, let $U' = P \times^{S_d} \{ 1,
\ldots, d\}$. Then $U' \to U$ is \'{e}tale of degree $d$. One can recover $P$
from $U'$ as the sheaf on the small \'{e}tale site of $U$,
\begin{equation*}
  P = \op{Iso}_{et}(U', U \times \{ 1, \ldots, d\})
\end{equation*}
Observe that $P$ is an \'{e}tale locally trivial principal $S_d$ bundle. This
is because the map $U' \to U$ is \'{e}tale locally isomorphic to $U \times \{1
\ldots d \} $ - proper, separated, surjective, \'{e}tale maps of degree $d$ are
precisely the maps with this property.  Then,
\begin{equation*}
  \op{Hom}(P, X^d )^{S_d}= \op{Hom}(P\times\{1, \ldots , d\}, X )^{S_d} =
  \op{Hom}(P\times^{S_d} \{ 1, \ldots , d\} , X) = \op{Hom}(U', X)
\end{equation*}
\end{proof}

\begin{corollary}
Let $V$ be a DM stack.  The $2$-groupoid $\Hom(V, S^d X)$ is equivalent to the
$2$-groupoid of diagrams $V \from V' \to X$, with $V' \to V$ proper, separated
surjective, and \'{e}tale of degree $d$. Further, the $2$-groupoid
$\op{HomRep}(V, S^d X)$ of representable maps $V \to S^d X$ is equivalent to
the $2$-groupoid of such diagrams $V \from V' \to X$, where the inertia groups
of $V$ act faithfully on the fibres of $V' \to V$.
\end{corollary}
\begin{proof}
We prove the statement about representability.  $V \to S^d X$ is representable
if and only if the principal $S_d$-bundle, $P \to V$, is an algebraic space.
This is equivalent to saying that the inertia groups of $V$ act faithfully on
the fibres of $V' \to V$.
\end{proof}

\subsection{Stacks of \'{e}tale covers}
We need some notation to shorten the cumbersome $g,I,m,a$ labels.  Let
$\vertex(A)$ be the groupoid of quadruples $\nu = (g(\nu), T(\nu), m, a(\nu) )$
where $g(\nu) \in \Z_{\ge 0}$, $T(\nu)$ is a finite set, $m: T(\nu) \to \Z_{>
0}$ is a function, and $a(\nu) \in A$. We impose the stability condition as
before: if $a(\nu) = 0$, then $2g(\nu) - 2 + \#T(\nu) > 0$. The morphisms are
isomorphisms preserving all the structure. Let $\rvertex(A)$ be the groupoid of
triples $v = (g(v), T(v), a(v))$ satisfying the stability condition. There is a
map $r: \vertex(A) \to \rvertex(A)$ sending $(g,I,m,a) \to (g,I,a)$.   For $\nu
\in \vertex$ we have the moduli stack $\mod_\nu$; similarly for $v \in
\rvertex$ we have $\mod_v$. There is  a map $r: \mod_\nu \to \mod_{r(\nu)}$
which associates to a twisted curve its coarse moduli space.

We also want labels for moduli stacks of disconnected curves.  We define a
groupoid $\op{exp}(\vertex)$. An object $\alpha \in \op{exp}(\vertex)$, is  a
finite set $V(\alpha)$, and a map $V(\alpha) \to \op{Ob} \vertex$. A morphism
$\alpha \to \alpha'$ in $\op{exp} (\vertex)$, is an  isomorphism $\phi :
V(\alpha) \iso V(\alpha')$ of finite sets, together with an isomorphism $v \iso
\phi(v)$ of the associated element of $\vertex$, for each $v \in V(\alpha)$.
Define $\op{exp} (\rvertex)$ in a similar way. Given $\alpha \in \op{exp}
(\vertex)$, let
\begin{equation*}
  T(\alpha) = \coprod_{v \in V(\alpha)} T(v)
\end{equation*}
The groupoid $\op{exp} (\vertex)$ labels possibly disconnected nodal curves.
Let
\begin{equation*}
  \mod_\alpha = \prod_{v \in V(\alpha)} \mod_{v}
\end{equation*}

Next, we want to define labels for \'{e}tale maps of twisted curves, $\mc C'
\to \mc C$. $\mc C'$ may be disconnected. A covering $\eta$, consists of
\begin{enumerate}
  \item An element $s(\eta) \in \op{exp} (\vertex)$, the source, and an element $t(\eta) \in
  \vertex$, the target.
  \item   These must satisfy

\begin{equation*}
  \sum_{v \in V(s(\eta))} a(v) = a(t(\eta)) \in A
\end{equation*}
    \item A map of finite sets, $p: T(s(\eta)) \to T(t(\eta))$.
  \item For each $t' \in T(s(\eta))$, we require that $m(t')$ divides $m(p(t'))$. Let $d(t') =
  m(p(t'))/m(t')$.
  \item We require that for each $t \in T(t(\eta))$,
\begin{equation*}
  m(t) = \op{lcm}\{ d(t')
  \mid p(t') = t \}
\end{equation*}
   where $\op{lcm}$ stands for lowest common multiple.
  \item A function $d: V(s(\eta)) \to \Z_{\ge 1}$, the degree.
  \item For each $t \in T(t(\eta))$, and each $v \in V(s(\eta))$,
\begin{equation*}
  \sum_{\substack{ t' \in T(v)\\ p(t') = t}} d(t') = d(v)
\end{equation*}
  We define $d(\eta) = \sum_{v \in V(s(\eta))} d(v)$.
  \item The Riemann-Hurwitz formula holds: for each $v \in V(s(\eta))$,
\begin{equation*}
  2(g(s(\eta)_v) -1) = 2d(v)(g(t(\eta)) -1) + \sum_{t' \in T(v)} ( d(t') - 1)
\end{equation*}
\end{enumerate}
 Let $\cvertex$ be the groupoid of all coverings $\eta$, with the obvious
isomorphisms.   There are source and target functors,
\begin{align*}
  s &: \cvertex \to \op{exp} (\vertex) \\
  t &: \cvertex \to \vertex
\end{align*}
We will often write $\alpha \cover \nu$ to mean a covering $\eta$ with $s(\eta)
= \alpha$ and $t(\eta) = \nu$.

Associated to a covering $\eta \in \cvertex$, we define a stack $\mod_{\eta}$
of \'{e}tale covers $f: \mc C' \to \mc C$. $\mod_\eta$ is the category whose
objects are
\begin{itemize}
  \item An object of $\mod_{t(\eta)}$, with associated family of twisted nodal curves
  $\mc C \to U$, sections $T(t(\eta)) \to \mc C$ and constructible function $f:
  C_{\op{gen}} \to U$, where $C$ is the coarse moduli space of $\mc C$.
  \item
  An object of $\mod_{s(\eta)}$, with associated family of possibly disconnected
  twisted nodal curves
  $\mc C' \to U$, sections $T(s(\eta)) \to \mc C'$ and constructible function $f':
  {C'}_{\op{gen}} \to U$, where $C'$ is the coarse moduli space of $\mc C$.
  \item
  An etale map $p: \mc C' \to \mc C$.
\end{itemize}
These must satisfy:
\begin{itemize}
\item
The diagram
\begin{equation*}
\begin{diagram}
  \node{T(s(\eta)) \times U} \arrow{e} \arrow{s} \node{\mc C'} \arrow{s} \\
  \node{T(t(\eta)) \times U} \arrow{e} \node{\mc C}
\end{diagram}
\end{equation*}
must be Cartesian over $U$. This implies, in particular, that the marked points
of $\mc C'$ are precisely those lying over marked points of $\mc C$.
\item
Let $p_\ast f'$ be the constructible function on $C_{\op{gen}}$ given by
pushing forward $f'$; we require that $p_\ast f' = f$.
\item
The map $\mc C \to BS_d$ associated to the etale map $\mc C' \to \mc C$ must be
representable.
\end{itemize}
\begin{proposition}
\label{prop mod algebraic}
The map $\mod_\eta \to \mod_{t(\eta)}$ is \'{e}tale; therefore, $\mod_\eta$ is
smooth. Further, $\mod_\eta$ is algebraic.
\end{proposition}
Before we prove this, we need a lemma.

\begin{lemma}
Let $G$ be a finite group.  The stack
\begin{equation*}
  \mod_{g,I,m}(BG) \defeq \op{HomRep}_{\mod_{g,I,m}}(\curve_{g,I,m}, BG \times
  \mod_{g,I,m})
\end{equation*}
of representable maps $\mc C \to BG$ from curves $\mc C \in \mod_{g,I,m}$ is
algebraic.
\end{lemma}
\begin{proof}[Sketch of proof]
A representable map from a twisted nodal curve $\mc C$ to $BG$ is the same as a
principal $G$ bundle $P \to \mc C$, whose total space is an ordinary nodal
curve.  This is the same, just as in \cite{acv}, Theorem 4.3.2, as a nodal
curve $P$, with a $G$-action, such that the map $P \to P/G$ to the scheme
quotient is generically a principal $G$-bundle;  the $G$ action must also have
some compatibility at the nodes. We recover $\mc C$ as $[P/G]$, the stack
quotient.

The stack of nodal curves $P$ is algebraic. Further, the stack of nodal curves
with a $G$ action is algebraic, because a $G$ action on a curve $P$ can be
identified with its graph in $P \times P \times G$.  It follows that
$\mod_{g,I,m}(BG)$ is algebraic.
\end{proof}
\begin{proof}[Proof of Proposition \ref{prop mod algebraic}]
The deformations of an \'{e}tale cover $\mc C' \to \mc C$ are the same as
deformations of the base, as in \cite{acv}; therefore the map $\mod_\eta \to
\mod_{t(\eta)}$ is \'{e}tale.

 Consider the map $\mod_\eta \to
\mod_{g(t(\eta)), T(t(\eta)), m}(BS_d)$.  This map is relatively a scheme; it
follows that $\mod_\eta$ is algebraic.
\end{proof}
The fact that $\mod_\eta$ is algebraic implies that its image in
$\mod_{t(\eta)}$, which is an open substack, is also algebraic.

Let us look at coverings $\mc C' \to \mc C \in \mod_\eta$ locally.  Given  a
tail $t \in T(t(\eta))$, in an etale neighbourhood of the twisted marked point
$t \to \mc C$, $\mc C' \to \mc C$ looks like
\begin{equation*}
 \left( \Spec \left(\oplus_{t' \in p^{-1}(t)} \oplus_{i = 1}^{d(t')} k[x_{t',i}]\right) \to \Spec k[y] \right) / \mu_{m(t)}
\end{equation*}
As before, $d(t') = m(t)/m(t')$. The algebra map which induces this map of
schemes is of course $y \to \sum x_{t',i}$. The $\mu_{m(t)}$ action sends $y
\to ly$, for $l \in \mu_{m(t)}$; and $x_{t',i} \to lx_{t',i+1 \op{mod} d(t')
}$. This action is faithful, which is equivalent to representability of the
associated map $\mc C \to BS_d$,  because $m(t)$ is the lowest common multiple
of $d(t')$ for $t' \in p^{-1}(t)$. The stabilizer of $x_{t',i}$ is $\mu_{m(t')}
\subset \mu_{m(t)}$. There is a similar picture near the nodes, except $k[x]$
is replaced by $k[u,v]/uv$.

\subsection{Stacks of stable maps}
Let $X$ be a smooth projective variety.  We will let our semigroup $A$, be the
Mori cone of effective $1$-cycles on $X$, up to numerical or homological
equivalence.  For each $v = (g,I,\beta) \in \rvertex$, we have the stack
$\cmod_v(X)$ of Kontsevich stable maps in $X$.  There is a map
\begin{equation*}
  \cmod_v(X) \to \mod_v
\end{equation*}
which associates to a stable map $C \to X$ with marked points, the curve $C$,
with its marked points, and the constructible function $C_{\op{gen}} \to C(X)$
given by taking the homology class of an irreducible component.  In a similar
way, for each $\alpha \in \op{exp}(\rvertex)$, labelling disconnected curves,
we have a moduli stack $\cmod_\alpha(X)$ with a map $\cmod_\alpha(X) \to
\mod_\alpha$.
\begin{lemma}
\begin{equation*}
  \cmod_{g,n,\beta}(X) =  \mod_{g,n,\beta} \times_{\mod_{g,n-1,\beta}}
\cmod_{g,n-1,\beta}(X)
\end{equation*}
\end{lemma}
\begin{proof}
It was shown in \cite{behrend manin}  that $\cmod_{g,n,\beta}(X)$ is the
universal curve over $\cmod_{g,n-1,\beta}(X)$, which is pulled back from
$\mod_{g,n-1,\beta}$. But we have shown that $\mod_{g,n,\beta}$ is the
universal curve over $\mod_{g,n-1,\beta}$.
\end{proof}
More generally,
\begin{equation*}
  \cmod_{g,n,\beta}(X) = \mod_{g,n,\beta} \times_{\mod_{g,\beta}} \cmod_{g,\beta}(X)
\end{equation*}
so that all stacks of marked stable maps to $X$, arise by base change with the
stack of unmarked stable maps $\cmod_{g,\beta}(X)$.

Let $V$ be a proper projective Deligne-Mumford stack.  Abramovich and Vistoli
\cite{av} defined the stack of stable maps to $V$: this is the stack of
representable maps $f: \mc C \to V$ from twisted nodal curves with marked
points to $V$, such that $\op{Aut}(f)$ is finite.  We are only interested in
the case $V = S^d X$, the symmetric product stack of a smooth projective
variety $X$.  Take our semigroup to be $C(X)$ as above.  For each $\nu  =
(g,I,m,\beta) \in \vertex$, let $\cmod_\nu(S^d X)$ be the stack of stable maps
from curves $\mc C \in \mod_{g, I, m}$, such that if $\mc C \from \mc C' \to X$
is the associated correspondence, then $\mc C' \to X$ has class $\beta \in
C(X)$.

For each covering $\eta \in \cvertex$, define
\begin{equation*}
  \cmod_{\eta}(X) = \mod_\eta \times_{\mod_{s(\eta)}} \cmod_{s(\eta)}
\end{equation*}
Let $\op{Aut}(\eta \mid t(\eta))$ be the group of automorphisms of $\eta$ which
act trivially on $t(\eta)$.
\begin{lemma}
There is a natural isomorphism
\begin{equation*}
 \coprod_{\eta, \, t(\eta) = v} \cmod_{\eta}(X) / \op{Aut}(\eta \mid t(\eta)) \iso   \cmod_v(S^d X)
\end{equation*}
\end{lemma}

\subsection{Tautological line bundles}
Let $\nu \in \rvertex$ or $\vertex$.  For each $t \in T(\nu)$, there is a
section $\sigma_t: \mod_\nu \to \curve_\nu$ of the universal curve.  Define
$L_t = \Omega^1_{\sigma_t}$ to be the relative cotangent bundle.  $L_t$ is the
tautological line bundle. If $\mu \in \vertex$, so that $r(\mu) \in \rvertex$,
we have a map $\mod_\mu \to \mod_{r(\mu)}$.  For each $t \in T(\mu) =
T(r(\mu))$, we have $r^\ast L_t = L_t^{\otimes m(t)}$.

 For $\eta \in \cvertex$, for each $t \in T(s(\eta))$ (or $t \in T(t(\eta))$) there is a tautological
line bundle $L_t$, pulled back from $\mod_{s(\eta)}$ (respectively
$\mod_{t(\eta)}$). If $p(t') = t$ under the projection $T(s(\eta)) \to
T(t(\eta))$, then $L_{t'} \iso L_t$.

Let $\psi_t = c_1(L_t)$ on any of the three types of moduli stack.

\subsection{Automorphisms and deformations of twisted nodal curves}
\label{automorphisms and deformations}
Let $v \in \vertex$, let $\mc C \in \mod_v$ and let $C \in \mod_{r(v)}$ be the
coarse moduli space of $\mc C$.  We want to describe $\op{Aut}(\mc C \mid C)$,
the group of automorphisms of $\mc C$ which are trivial on the coarse moduli
space $C$. This group splits as a product of contributions from each twisted
node and twisted marked point of $\mc C$:  a twisted node or marked point with
inertia group $\mu_r$ contributes $\mu_r$.

 The fibre of $\mod_v \to \mod_{r(v)}$ over a curve $C \in
\mod_{r(v)}$ can similarly be described as a product of local contributions
from the nodes and marked points of $C$.  For each tail $t \in T(v) = T(r(v))$,
we have a factor of $B \mu_{m(t)}$.  For each node of $C$, we have a factor of
$\coprod_{k \in \Z_{> 0}} B\mu_k$.   For more details, see \cite{acv}.

For $v \in \vertex$, let $T_1(v) \subset T(v)$ be the set of tails with
multiplicity $m(t) = 1$.  For $\mc C \in \mod_v$ and $t \in T_1(v)$, the marked
point $P_t \in \mc C$ is untwisted.  The first-order deformations of $\mc C$
are given by
\begin{equation*}
  \op{Ext}^1(\Omega^1_{\mc C}(\sum_{t \in T_1(v)}P_t), \Oo_{\mc C})
\end{equation*}
The deformation theory is unobstructed.  We can identify this space with
$$H^0 \left( \omega_{\mc C} \otimes \Omega^1_{\mc C}( \sum_{t \in T_1(v)} P_t)
\right)^\vee$$
where $\omega_{\mc C}$ is the dualizing line bundle.

Let $C \in \mod_{r(v)}$ be the coarse moduli space of $\mc C$. We have a map
$\pi : \mc C \to C$.  In \cite{av}, it is shown that $\pi_\ast$ is an exact
functor, and so
\begin{equation*}
 H^0 \left(\omega_{\mc C} \otimes \Omega^1_{\mc C}\left( \sum_{t \in T_1(v)}
 P_t\right)\right)
 = H^0 \left(\pi_\ast \left( \omega_{\mc C} \otimes \Omega^1_{\mc C}\left( \sum_{t \in
 T_1(v)} P_t\right)\right)\right)
\end{equation*}
The space of first order deformations of $C$ is similarly given by $H^0
\left(\omega_{C} \otimes \Omega^1_{C}( \sum_{t \in T(r(v))}
 P_t)\right)^\vee$.  Observe that we have a pole at all tails, not just those with
multiplicity one.   Clearly
\begin{equation*}
  \pi_\ast \left( \omega_{\mc C} \otimes \Omega^1_{\mc C}\left( \sum_{t \in
 T_1(v)} P_t\right)\right) = \omega_C \otimes \Omega^1_C\left(\sum_{t \in T(r(v))} P_t \right)
\end{equation*}
away from the nodes and marked points of $C$.  In fact, this equality extends
also to the marked points. At the nodes, however, this is no longer true. The
map $\mod_v \to \mod_{r(v)}$ is ramified along the divisor of singular curves.
The degree of ramification along the divisor along the divisor in $\mod_v$
corresponding to a node with inertia group $\mu_k$ is $k-1$ (i.e.\ a function
vanishing to degree $1$ on the divisor in $\mod_{r(v)}$ vanishes to degree $k$
along the divisor in $\mod_v$ when its pulled back).  One can see this by
looking at the local picture, as in \cite{acv}, section 3.

Let $\eta \in \cvertex$ and let $\mc C' \to \mc C \in \mod_\eta$ be an
\'{e}tale cover of twisted balanced curves.  Let $C' \to C$  be the
corresponding ramified covering of the coarse moduli spaces of $\mc C'$ and
$\mc C$. We are interested in $\op{Aut}(\mc C' \to \mc C \mid C' \to C)$, the
automorphisms which are trivial on the coarse moduli space.  As before, this
splits as a product with a contribution from each twisted node and twisted
marking. Each $t \in T(t(\eta))$ -- that is each marking of $\mc C$ --
contributes $\mu_{m(t)}$. However, in this case the contribution from the nodes
is trivial. This follows from the fact that the the map $\mc C \to BS_d$ is
representable. The marked points of $\mc C'$ do not contribute anything extra.

We have seen already what the deformations of $\mc C' \to \mc C$ are: they are
the same as deformations of $\mc C$, that is the map $\mod_{\eta} \to
\mod_{t(\eta)}$ is \'{e}tale.

\section{Graphs}
\label{sec_graphs}
We define categories of graphs, which label various flavours of nodal curve, as
well as \'{e}tale covers of twisted nodal curves. We introduce three
categories: $\rgraph$, which contains labels for untwisted nodal curves;
$\graph$, which has labels for twisted nodal curves; and $\cov$, which has
labels for pairs of twisted nodal curves $\mc C'$, $\mc C$, with an \'{e}tale
map $\mc C' \to \mc C$. These categories depend on a semigroup $A$, which we
will usually not mention. The morphisms in these categories correspond to
degenerating curves, or dually to contracting graphs. There are functors,
denoted $\mod$, from each of these categories to the category of stacks, which
take a label to the moduli stack of all curves with that label; as well as
functors
\begin{equation*}
 \cov \overset{s,t}{\rightrightarrows} \graph \overset{r}{\to} \rgraph
\end{equation*}
$s$ and $t$ stand for source and target, and take the labels for a pair $\mc C'
\to \mc C$ to the labels for $\mc C'$ or $\mc C$ respectively.  $r$ takes the
labels for $\mc C$ to the labels for its coarse moduli space $C$.  These
functors get translated into morphisms of stacks after applying the moduli
stack functor $\mod$; that is, there are maps of stacks $s: \mod_{\sigma} \to
\mod_{s(\sigma)}$, and similarly for $t$ and $r$.

Now we define these categories. Let $\graph$ be the category whose objects are
objects  $\eta$ of  $\op{exp} (\vertex)$, together with an order two
isomorphism $\sigma : T(\eta) \iso T(\eta)$, commuting with the multiplicity
map $m: T(\eta) \to \Z_{\ge 1}$. Define the objects of $\rgraph$ in a similar
fashion, using $\rvertex$ instead of $\vertex$ and omitting references to $m$.
There are forgetful maps $F: \graph \to \op{exp}(\vertex)$, and $F: \rgraph \to
\op{exp}(\rvertex)$.  For $\gamma \in \graph$ (or $\gamma \in \rgraph$), the
vertices of $\gamma$, written $V(\gamma)$, is the set of vertices
$V(F(\gamma))$ of the underlying element of $\op{exp}(\vertex)$. The half-edges
of $\gamma$ is the set $H(\gamma) = T(F(\gamma))$. The set of edges of
$\gamma$, $E(\gamma)$, is the set of free $\Z/2$ orbits on $H(\gamma)$.  The
set of tails of $\gamma$, $T(\gamma)$, is the set of $\sigma$-fixed points on
$H(\gamma)$. We can fit the structure of a graph $\gamma \in \graph$ into the
diagram
\begin{equation*}
  \begin{diagram}
  \node{\Z_{> 0}}  \node{H(\gamma)\circlearrowleft^\sigma} \arrow{w,t}{m}\arrow{s}  \\ 
  \node{\Z_{\ge 0}}  \node{V(\gamma)} \arrow{w,t}{g} \arrow{e,t}{a} \node{A}
  \end{diagram}
\end{equation*}

Given $\gamma \in \graph$ or $\rgraph$, let $\mod_\gamma = \prod_{v \in
V(\gamma)} \mod_v$.  For every half-edge $h \in H(\gamma)$, there is a
tautological line bundle $L_h$.  For every edge $e \in E(\gamma)$ corresponding
to the $\sigma$-orbit $(h_1,h_2)$ let $L_e = L_{h_1} \otimes L_{h_2}$.

\subsection{Contractions}
Now we want to define the morphisms in the categories $\graph$ and $\rgraph$,
called \emph{contractions}. In terms of the cell complex $C(\gamma)$ associated
to a graph $\gamma$, a contraction $\gamma' \cont \gamma$ is a surjective
continuous map $C(\gamma') \to C(\gamma)$ which is an isomorphism away from the
vertices of $C(\gamma)$, and possibly maps some edges of $C(\gamma')$ to
vertices of $C(\gamma)$.    It is better to describe contractions more
formally. Let $\gamma, \gamma' \in \graph$. A contraction $\gamma' \cont
\gamma$ is a surjective map of sets $f: H(\gamma') \coprod V(\gamma') \to
H(\gamma) \coprod V(\gamma)$, such that
\begin{itemize}
\item
    $V(\gamma') \subset f^{-1}(V(\gamma))$.
\item
    The diagram
\begin{equation*}
  \begin{diagram}
    \node{H(\gamma') \coprod V(\gamma')} \arrow{e,t}{f} \arrow{s,l}{\sigma \coprod 1} \node{ H(\gamma) \coprod
V(\gamma)}  \arrow{s,r}{\sigma \coprod 1} \\
    \node{H(\gamma') \coprod V(\gamma')} \arrow{e,t}{f} \arrow{s} \node{ H(\gamma) \coprod
V(\gamma)}  \arrow{s} \\
    \node{V(\gamma')}  \arrow{e,t}{f} \node{V(\gamma)}
  \end{diagram}
\end{equation*}
    commutes, where $\sigma$ is the involution on $H(\gamma')$ or $H(\gamma)$,
    and $H(\gamma) \coprod V(\gamma) \to V(\gamma)$ comes from the map $H(\gamma)
    \to V(\gamma)$ assigning to a half-edge the vertex it is attached to.
\item
    $f$ induces an isomorphism $H(\gamma') \supset f^{-1}(H(\gamma)) \iso
    H(\gamma)$, commuting with the multiplicity functions.
\item
    We require that $f$ does not contract any tails to vertices, so that $T(\gamma') =
    H(\gamma')^\sigma \subset f^{-1}(H(\gamma))$.  This implies that $f$
    induces an isomorphism $T(\gamma') \iso T(\gamma)$.
\item
    For each $v \in V(\gamma)$, we can define a graph $\gamma'_v \in \graph$,
    with vertices $f^{-1}(v) \cap V(\gamma')$, edges $f^{-1}(v) \cap
    E(\gamma')$, and tails $f^{-1}(H(v)) \coprod (f^{-1}(v) \cap T(\gamma'))$.
    We require that $\gamma'_v$ is connected of genus $g(\gamma'_v) = g(v)$.
\end{itemize}
We define contractions of graphs $\gamma, \gamma' \in \rgraph$ in the same
fashion, leaving out references to the multiplicity function.

Let $E(f) \subset E(\gamma')$ be the set of edges contracted by $f$.  Given any
subset $I \subset E(\gamma')$ there is a unique contraction $f: \gamma' \to
\gamma$ with $E(f) = I$.  One can think of $E(f)$ as being the kernel of $f$.

Contractions correspond to degenerating curves by adding more nodes.  If we
have a contraction $f: \gamma' \cont \gamma$, we can identify $\mod_{\gamma'} =
\prod_{v \in V(\gamma)} \mod_{\gamma'_v}$.  There are maps $\mod_{\gamma'_v}
\to \mod_{v}$, which induces a map
$$f_\ast : \mod_{\gamma'} \to \mod_{\gamma}$$
The map $f_\ast$ has cotangent complex on $\mod_{\gamma'}$
\begin{equation*}
f^\ast \Omega^1_{\mod_{\gamma}} \to   \Omega^1_{\mod_{\gamma'}}
\end{equation*}
We can compute the cohomology of this complex in terms of tautological line
bundles on $\mod(\gamma')$.
\begin{align*}
\op{Ker} \left(   f^\ast \Omega^1_{\mod_{\gamma}} \to
    \Omega^1_{\mod_{\gamma'}} \right) &= \oplus_{i \in E(f)} L_i \\
\op{CoKer}\left(   f^\ast \Omega^1_{\mod_{\gamma}} \to
    \Omega^1_{\mod_{\gamma'}} \right) &= 0
\end{align*}

\subsection{Coverings}
\label{subsection coverings}
We define labels for \'{e}tale maps $\mc C' \to \mc C$ of twisted nodal curves.
A covering $\gamma' \cover \gamma$ is a map of sets $p : V(\gamma') \to
V(\gamma)$, and for each $v \in V(\gamma)$, a covering $p^{-1}(v) \to v$. Here
we consider $p^{-1}(v)$ as an object of $\op{exp}(\vertex)$. We require that
the associated map $H(\gamma') \to H(\gamma)$ is equivariant with respect to
the involution, takes tails to tails, and edges to edges.

Given a covering $p:\gamma' \cover \gamma$, for $v \in V(\gamma)$, let
$\gamma'_v = p^{-1}(v) \in \op{exp}(\vertex)$ be the vertices lying over $v$.
We define a stack $\mod_{\gamma' \cover \gamma}$ of coverings, by
\begin{equation*}
  \mod_{\gamma' \cover \gamma} = \prod_{v \in V(\gamma)} \mod_{\gamma'_v \cover v}
\end{equation*}
As usual there are maps $\mod_{\gamma'\cover\gamma} \to \mod_{\gamma}$, which
is \'{e}tale, and $\mod_{\gamma'\cover \gamma} \to \mod_{\gamma'}$.  If we also
have a contraction $\gamma \cont \eta$, there is a unique covering $\eta'
\cover \eta$ and a contraction $\gamma' \cont \gamma$ such that the diagram
\begin{equation*}
  \begin{diagram}
    \node{\gamma'} \arrow{e} \arrow{s,A} \node{\eta'} \arrow{s,A} \\
    \node{\gamma} \arrow{e} \node{\eta}
  \end{diagram}
\end{equation*}
commutes.

Let $\cov$ be the category whose objects are coverings $\gamma' \cover \gamma$
with $\gamma', \gamma \in \graph$ and whose morphisms are diagrams as above.
There are source and target functors $s,t: \cov \to \graph$. There is a functor
$\cov \to \text{stacks}$, sending a covering $\rho = s(\rho)\cover t(\rho)$ to
$\mod_\rho$, and an arrow $\rho \cont \phi$ to the associated map of stacks
$\mod_\rho \to \mod_\phi$.

Let $\rho,\phi \in \cov$ and suppose we have a morphism $f: \rho \to \phi$.
There is an induced morphism $t(\rho) \to t(\phi)$: let $E(f) \subset
E(t(\rho))$ be the set of edges contracted.  There is a map $\mod_\rho \to
\mod_\phi$, which is \'{e}tale over its image.  The relative cotangent bundle
is $\Omega^1_{f_\ast} = \oplus_{i \in E(f)} L_i$.  The diagram
\begin{equation*}
  \begin{diagram}
  \node{\mod_\rho} \arrow{e}\arrow{s} \node{\mod_\phi}
  \arrow{s} \\
  \node{\mod_{t(\rho)}} \arrow{e} \node{\mod_{t(\phi)}}
  \end{diagram}
\end{equation*}
commutes, and the vertical maps are \'{e}tale.

\section{Pull backs of tautological classes}
\label{sec_pull_back}

We want to do intersection theory on our moduli stacks, later with virtual
fundamental classes. We will work both with Artin algebraic stacks and with
Deligne-Mumford stacks. For DM stacks, we will use the theory $A_\ast$ of
Vistoli \cite{vistoli}, with $\Q$-coefficients.  For Artin stacks, Kresch has
defined in \cite{kresch} intersection theory, as has Toen in \cite{toen riemann
roch}.  However, all we ever need to do with Artin stacks is intersect
regularly embedded smooth divisors with regularly embedded smooth closed
substacks, and also take first Chern classes of line bundles. No particularly
sophisticated technology is required -- what we need is basically the same as
that used by Behrend in \cite{behrend gw} to define Gromov-Witten classes.

We want to intersect cycles corresponding to graphs. If we have a map $f: \rho'
\to \rho$ in one of our categories of graphs $\cov,\graph,\rgraph$ , we have an
associated map $f_\ast : \mod_{\rho'} \to \mod_{\rho}$. Let $\mod_f \into
\mod_{\rho}$ be the closed substack of $\mod_{\rho}$ which is supported on the
image of $f$.  $\mod_f$ is smooth, and the map $\mod_f \into \mod_\rho$ is a
closed regular embedding. The map $\mod_{\rho'} \to \mod_f$ is \'{e}tale, in
general non-representable, of degree $\norm f$, where:
\begin{enumerate}
  \item If $f: \rho' \to \rho$ is a morphism in $\rgraph$, then
\begin{equation*}
    \norm {f} = \#\Aut(\rho' \to \rho \mid \rho)
\end{equation*}
 is the number of automorphisms of
$\rho'$ commuting with the map $\rho' \to \rho$.
  \item If $f: \rho' \to \rho$ is a morphism in $\graph$, then
\begin{equation*}
    \norm {f} = \frac{\#\Aut(\rho' \to \rho \mid \rho)} {\prod_{i \in E(f)} m(i)}
\end{equation*}
    where $E(f) \subset E(\rho')$ is the set of edges contracted.
  \item If $f: \rho' \to \rho$ is a morphism in $\cov$, then
\begin{equation*}
  \norm {f} = \frac{\#\Aut(\rho' \to \rho \mid \rho)} {\prod_{i \in E(f)} m(i)^2}
\end{equation*}
 where $E(f) \subset E(t(\rho'))$ is the set of edges contracted.
\end{enumerate}
Let $e \in E(f)$. Observe that the line bundle $L_e$ on $\mod_{\rho'}$ descends
to a line bundle $L_e$ on $\mod_f$.  Let $\psi_e = c_1(L_e) \in A^1 \mod_f$.
The conormal bundle to the embedding $\mod_f \into \mod_{\rho}$ is $\oplus_{e
\in E(f)} L_e$.

Let $\mc A$ be one of the categories $\cov, \graph,\rgraph$. Suppose we have a
diagram
\begin{equation*}
\begin{diagram}
    \node[2]{\rho''} \arrow{s,r}{g} \\
    \node{\rho'} \arrow{e,t}{f} \node{\rho}
\end{diagram}
\end{equation*}
in $\mc A$. We want to calculate $f^\ast [\mod_{g}] \in A_\ast \mod_{\rho'}$.
For simplicity we will assume that $\mod_g \into \mod_\rho$ is a divisor, or
equivalently $\# E(g) = 1$.  This the only case we will need.
\begin{lemma}
\label{lemma intersect graph}
$f^\ast [\mod_{g}]$ can be expressed as a sum with a term for each map $h :
\rho' \to \rho''$ such that $g \circ h = f$. Each such term of this is weighted
by $-\psi_e$, where $e$ is the unique element of $E(f) \setminus E(h)$. There
is also a term $[\mod_h]$ for each  isomorphism class of commutative diagrams
\begin{equation*}
  \begin{diagram}
    \node{\eta} \arrow{e,t}{k}\arrow{s,l}{h} \node{\rho''} \arrow{s,r}{g} \\
    \node{\rho'} \arrow{e,b}{f} \node{\rho}
  \end{diagram}
\end{equation*}
where $\#E(h) = 1$ and $E(h) \not\subset E(k)$.
\end{lemma}
The factor $-\psi_e$ is the first Chern class of the normal bundle. Now we
calculate pull backs under the maps of stacks induced by the functor $r: \graph
\to \rgraph$.
\begin{lemma}
Let $\gamma \in \graph$, and suppose $f: \alpha \to r(\gamma)$ is a morphism in
$\rgraph$, with $\# E(f) = 1$. Let $r: \mod_\gamma \to \mod_{r(\gamma)}$ be the
canonical map.
\begin{equation*}
  r^\ast [\mod_{f}] = \sum_{g: \gamma' \to \gamma} m(e) [\mod_{g}]
\end{equation*}
where the sum is over $g : \gamma' \to \gamma$ such that $r(g) = f$, and $e \in
E(g)$ is the unique element.
\end{lemma}
The factors $m(e)$ come from the fact that the map $\mod_\gamma \to
\mod_{r(\gamma)}$ has ramification along the boundary divisors of
$\mod_{r(\gamma)}$.

Next we calculate pullbacks under $s$. Let $\rho \in \cov$, and let $f : \gamma
\to s(\rho)$ be a morphism in $\graph$, where again $\# E(f) = 1$.  There is a
map $s: \mod_{\rho} \to \mod_{s(\rho)}$.
\begin{lemma}
$s^\ast [\mod_{f}] \in A_\ast \mod_\rho$  can be expressed as a sum over
isomorphism classes of pairs of maps $g: \rho' \to \rho$ in $\cov$, and $h:
s(\rho') \to \gamma$ in $\graph$, such that $f \circ h = s(g)$, and $\# E(g) =
1$:
\begin{equation*}
  s^\ast  [\mod_{f}] =
  \sum_{g,h} [\mod_g]
\end{equation*}
\end{lemma}
Observe that since $t: \mod_\eta \to \mod_{t(\eta)}$ is \'{e}tale, it is easy
to calculate pull backs under $t$.

\subsection{Pull backs of $\psi$-classes}
\label{pull back psi classes}
So far we have seen how to pull back divisors corresponding to graphs.  Next,
we want to pull back Chern classes of tautological line bundles.

Under a morphism $f: \gamma' \to \gamma$ in any of our categories $\cov,
\graph, \rgraph$, for any half-edge $h \in H(\gamma)$ there is a unique
half-edge $h' \in H(\gamma')$ with $f(h') = h$.  Then, $f^\ast \psi_h =
\psi_{h'}$.

For the functors $s,t: \cov \to \graph$, for any $\eta \in \cov$ and half edge
$h \in H(s(\eta))$, by definition $s^\ast \psi_h = \psi_h$; and similarly for
$t$. Also, if $h_s \in H(s(\eta))$ lies above $h_t \in H(t(\eta))$, then
$\psi_{h_s} = \psi_{h_t}$.

For the functor $r: \graph \to \rgraph$: if $\gamma \in \graph$, then
$H(\gamma) = H(r(\gamma))$. If $h \in H(\gamma)$, then $r^\ast \psi_h = m(h)
\psi_h$.

There is another case which is not so trivial.   Let $\eta \in \cvertex$.  Let
$I \subsetneq T(s(\eta))$ be such that after removing the tails in $I$,
$s(\eta)$ remains stable.  This always happens if $g(s(\eta))> 0$ or if
$a(s(\eta)) \in A$ is non-zero.   Let $v \in \rvertex$ be obtained from
$r(s(\eta))$ by removing the tails in $I$.  There is a map $\pi :
\mod_{r(s(\eta))} \to \mod_v$.  We are interested in calculating the pullbacks
of tautological $\psi$ classes in $\mod_v$, under the morphisms in the diagram
\begin{equation*}
\begin{diagram}
  \node{\mod_\eta} \arrow{e,t}{s}  \node{\mod_{s(\eta)}} \arrow{e,t}{r} \node{\mod_{r(s(\eta))}} \arrow{e,t}{\pi}
  \node{\mod_v}  \\
\end{diagram}
\end{equation*}
To do this we need to introduce yet more notation.

Let $t \in T(v)$ be any tail. Let $\gamma \in \rgraph$ with a contraction
$\gamma \to r(s(\eta))$. So, $I \subset T(\gamma)$ and $t \in T(\gamma)
\setminus I$. For $e \in E(\gamma)$, let $\gamma_e$ be the graph obtained by
contracting all edges except $e$.  We define $S(e,t,I) \in \{ 0, 1\}$ to be
$1$, if and only if $t$ is in a vertex of $\gamma_e$ that is contracted after
forgetting the tails $I$.  This happens if and only if $\gamma_e$ looks like
\begin{equation*}
  \begin{diagram}
        \\
     \node{\circ \, \, v_1}\arrow{w,b,-}{t} \arrow{e,b,-}{e} \arrow{n,l,-}{\text{tails }I' \subset I}
     \node{\circ \, \, v_2} \arrow{n,l,-}{\text{more tails}}
  \end{diagram}
\end{equation*}
where $v_1,v_2$ are the vertices of $\gamma_e$,  the genus $g(v_1) = 0$, and
the class $a(v_1) = 0 \in A$.  Define $S(e,t,I) = 0$ otherwise. Observe that
for each $e \in E(\gamma)$, there is at most one $t \in T(v) \into T(\gamma)$
such that $S(e,t,I) = 1$.

\begin{lemma}
For each $t \in T(v)$,
\begin{equation*}
  \pi^\ast \psi_t = \psi_t - \sum_{f: \gamma \to r(s(\eta))} [\mod_f]  S(e,t,I)
\end{equation*}
where the sum is over $f: \gamma \to r(s(\eta))$ with $\#E(\gamma) = \#E(f)=
1$, and $e \in E(\gamma)$ is the unique edge.
\end{lemma}
\begin{proof}
This is a rephrasing of a standard result.
\end{proof}
\begin{corollary}
For each $t \in T(v)$,
\begin{equation*}
 r^\ast \pi^\ast \psi_t = m(t) \psi_t - \sum_{f: \gamma \to s(\eta)} [\mod_f] m(e) S(e,t,I)
\end{equation*}
where the sum is over $f: \gamma \to s(\eta)$ with $\#E(\gamma) = \#E(f)= 1$,
and $e \in E(\gamma)$ is the unique edge.
\end{corollary}
\begin{corollary}
For each $t \in T(v)$,
\begin{equation*}
 s^\ast r^\ast \pi^\ast \psi_t = m(t) \psi_t -\sum_{f: \gamma \to \eta} [\mod_f] \sum_{e \in E(s(\gamma))} m(e) S(e,t,I)
\end{equation*}
where the sum is over $f: \gamma \to \eta$ with $\#E(t(\gamma)) = \#E(f)= 1$.
\label{cor s psi pullback}
\end{corollary}
For each tail $t \in T(v)$, let $z_t$ be a formal parameter. For each map
$\gamma \to \eta$ in $\cov$, and for each edge $e \in T(t(\eta))$ or tail $t
\in T(t(\eta))$, define formal parameters $w_e$ and $w_t$. We impose the
relations between $z$ and $w$ parameters:
\begin{itemize}
  \item For $t \in T(t(\eta))$, $w_t = \sum_{t' \in p^{-1}(t) \cap T(v)} m(t') z_{t'}$ where the sum is
  over tails of $T(v) \into T(s(\eta))$ lying over $t$.
  \item For $e \in E(t(\gamma))$, $w_e =  \sum_{e' \in p^{-1}(e)} \sum_{t \in T(v)} m(e') S(e',t,I) z_t$.
\end{itemize}
\begin{proposition}
With this notation,
\begin{equation}
s^\ast r^\ast \pi^\ast e^{\sum_{t \in T(v)} z_t \psi_t} = \sum_{f : \gamma \to
\eta}   i_\ast \left([\mod_f] e^{\sum_{t \in T(t(\gamma))} w_t \psi_t} \prod_{e
\in E(t(\gamma))} \frac{1 - e^{ w_e \psi_e} }{\psi_e}
   \right)
\label{eq psi pullback}
\end{equation}
where the sum is over isomorphism classes of maps $f: \gamma \to \eta$, and
$i: \mod_f \into \mod_\eta$ is the inclusion.
\end{proposition}
\begin{proof}
When all variables $z,w$ are zero, both sides are evidently equal.  Now apply
the operator $\frac{d}{dz_t}$ to both sides; it suffices to show that
$\frac{d}{dz_t}$ acts by intersection with $s^\ast r^\ast \pi^\ast \psi_t$ on
the right hand side.  This follows from corollary \ref{cor s psi pullback}, and
lemma \ref{lemma intersect graph}.
\end{proof}

\subsection{Pullbacks from Deligne-Mumford space}
\label{sec genus 0 relations} Let $v \in \rvertex$, be such that $g(v) = 0$ and $\#T(v) \ge 3$.  There
is a map
\begin{equation*}
  \pi : \mod_v \to \cmod_{0,T(v)}
\end{equation*}
where $\cmod_{0,T(v)}$ is the usual Deligne-Mumford stack of genus $0$ stable
curves.  This map is flat; this follows from the analogous result in
\cite{behrend gw}. Suppose $\#T(v) \ge 4$. For each distinct $i,j,k,l \in
T(v)$, there is a map
\begin{equation*}
  \cmod_{0,T(v)} \to \cmod_{0,\{i,j,k,l\}}
\end{equation*}
In the usual way, by pulling back two rationally equivalent divisors on
$\cmod_{0,\{i,j,k,l\}} \iso \mbb P^1$, we get the associativity equation on
$\mod_v$:
\begin{equation*}
  \sum_{f_{ij \mid jk} : \gamma \to v} [\mod_{f_{ij \mid kl}}] = \sum_{f_{ik \mid jl}: \gamma \to v}
  [\mod_{f_{ik\mid jl}}]
\end{equation*}
where the left hand side, is the sum over graphs $\gamma \to v$, such that
$\#E(\gamma) = 1$, and the tails $i,j$ and $k,l$ are on separate vertices of
$\gamma$; and similarly for the right hand side.

For each vertex $\nu \in \vertex$, with $g(\nu) = 0$ and $\#T(\nu) \ge 4$,
pulling this relation back from the map $\mod_\nu \to \mod_{r(\nu)}$, we get
the associativity relations on $\mod_{\nu}$:
\begin{equation*}
  \sum_{f_{ij \mid jk} : \gamma \to v} [\mod_{f_{ij \mid kl}}]m(e) = \sum_{f_{ik \mid jl}: \gamma \to v}
  [\mod_{f_{ik\mid jl}}]m(e)
\end{equation*}
where, as before, the graphs $\gamma_{ij \mid kl}$ have one edge, $e$, and the
tails $i,j$ and $k,l$ are in different vertices of $\gamma$.

For $v \in \rvertex$ with $g(v) = 0$ and $\#T(v) \ge 3$ as before, and for each
distinct $i,j,k \in T(v)$ consider the map $p: \cmod_{0,T(v)} \to
\cmod_{0,\{i,j,k\}}$. $\psi_i = 0$ on $\cmod_{0,\{i,j,k\}}$, because this is a
point. It follows that, on $\mod_{v}$, we have the equation
\begin{equation*}
  \psi_i = \sum_{f_{i \mid kl} : \gamma \to v} [\mod_{f_{i \mid kl}}]
\end{equation*}
where the sum is over graphs $f_{i \mid kl} : \gamma \to v$, such that $\# E(v)
= 1$, and the tails $i$ and $\{k,l\}$ are on different vertices.

Now let $\nu \in \vertex$, be such that $g(\eta) = 0$ and $\#T(\eta) \ge 3$.
Pulling back this relation from $\mod_{r(\nu)}$ to $\mod_{\nu}$,  gives us, for
each $i,j,k \in T(\nu)$,
\begin{equation*}
 m(i) \psi_i = \sum_{f_{i \mid kl} : \gamma \to v} [\mod_{f_{i \mid kl}}] m(e)
\end{equation*}
where the sum is over maps $f_{i, \mid kl} : \gamma \to v$ where $e \in
E(\gamma)$ is the unique edges, and the tails $i$ and $\{k,l\}$ are on
different vertices of $\gamma$.

Finally, observe that for each $\eta \in \cov$ with $t(\eta)$ having only one
vertex,  we can pull these relations back via the \'{e}tale map $\mod_\eta \to
\mod_{t(\eta)}$ in an obvious way.

\section{Virtual fundamental classes}
\label{sec virtual}
I will use the Behrend-Fantechi \cite{bf} construction of virtual fundamental
classes.

Let $F$ be a Deligne-Mumford stack, $V$ an Artin stack, and suppose there is a
map $F \to V$. A \emph{perfect relative obstruction theory} \cite{bf} is a
two-term complex of vector bundles $E = E^{-1} \to E^0$ on $F$, together with a
map
\begin{equation*}
  E \to \mbb L^\ast_{F/V}
\end{equation*}
in the derived category $D(F)$ to the relative cotangent complex $\mbb
L^\ast_{F/V}$, which is an isomorphism on $H^0$ and surjective on $H^{-1}$.

Associated to a perfect relative obstruction theory, Behrend-Fantechi in
\cite{bf} associate a cone $C \into   (E^{-1})^{\vee}$, and define the virtual
fundamental class
\begin{align*}
  [F,E] &\in A_\ast F \\
  [F,E] &= s^\ast [C]
\end{align*}
where $s: F \into E$ is the zero section. Here $s^\ast$ is the Gysin map, which
is defined to be the inverse of the pull-back isomorphism $\pi^\ast : A_\ast F
\to A_\ast E$.  Behrend-Fantechi show that $[F,E]$ only depends on the
quasi-isomorphism class of $E$.

Let me recall some of the details of their construction.  The \emph{relative
intrinsic normal cone} of a map $F \to V$, $\mf C_{F/V}$, is a cone stack over
$F$, with the property that if locally we factor $F \to V$ into $F
\overset{i}{\into} M \overset{p}{\to} V$, where $i: F \into M$ is a closed
embedding and $p: M \to V$ is smooth, then, $\mf C_{F/V}$ is the quotient stack
\begin{equation*}
  \mf C_{F/V} = [C_{F/M} / i^\ast T_{M/V}]
\end{equation*}
Here $C_{F/M}$ is the usual normal cone, which is acted on by the additive
group scheme $i^\ast T_{M/V}$.

Suppose we have a perfect relative obstruction theory $E \to \mbb
L^\ast_{F/V}$.  Let $E_1 = {E^{-1}}^\vee$ and $E_0 = {E^0}^\vee$.  One can show
that there is a closed embedding $\mf C_{F/V} \into [E_1/E_0]$, where
$[E_1/E_0]$ is the quotient stack. Form the Cartesian diagram
\begin{equation*}
\begin{diagram}
  \node{C} \arrow{e} \arrow{s} \node{E_1} \arrow{s} \\
  \node{\mf C_{F/V}} \arrow{e} \node{[E_1/E_0]}
\end{diagram}
\end{equation*}
where the vertical arrows are smooth, the horizontal arrows are closed
embeddings, and $C \to F$ is a usual cone, in particular a scheme over $F$. We
then define
\begin{equation*}
  [F,E] = s^\ast [C]
\end{equation*}
where $s: F \into E_1$ is the zero section.

Let $X', X$ be pure dimensional schemes of the same dimension, with $X$
irreducible, and let $f:X' \to X$ be a map between them. We say $f$ is of
degree $d$, if $f_\ast \Oo_{X'}$ is of rank $d$ over the generic point of $X$.
If $X$ is not irreducible, we say $f$ is of pure degree $d$ if it is of degree
$d$ for every irreducible component of $X$. This property is local in the
smooth topology of $X$.  That is, if $U \to X$ is a surjective smooth map, and
$U' = X' \times_X U$, then $U' \to U$ is of pure degree $d$ if and only if $X'
\to X$ is. Further, this property is local in the \'{e}tale topology of $X'$,
in the following sense. For each irreducible component $X_i$ of $X$, pick an
\'{e}tale cover $\coprod_j U_{ij} \to X'_i$, where $U_{ij}$ are connected and
$U_{ij} \to X'_i$ is of degree $e_{ij}$. . Then,
\begin{equation*}
  \op{deg}(X'_i/X_i) = \sum_j \frac{\op {deg} (U_{ij}/X_i)}{e_{ij}}
\end{equation*}
Let $V',V$ be Artin stacks of the same pure dimension, and let $V' \to V$ be a
map of relative Deligne-Mumford type. This means that  for every scheme $U \to
V$, $U \times_V V' \to U$ is a Deligne-Mumford stack.  We say $V' \to V$ is of
pure degree $d$, if for some smooth surjective map $U \to V$ from a scheme, for
each irreducible component $U_{i}$ of $U$, for some \'{e}tale atlas $\coprod_j
U'_{ij} \to U'_i = V' \times_V U_{i}$, with $U'_{ij} \to U'_i$ of degree
$e_{ij}$, $d = \sum_j \op {deg} (U'_{ij}/U_i)/e_{ij}$.

  This property
does not depend on the choices of smooth and \'{e}tale atlases, because for
schemes it is local in the smooth and \'{e}tale topologies, as above.  In
particular, if $V',V$ are DM stacks, then this definition, using the smooth
topology for $V$, agrees with the definition using the \'{e}tale topology for
$V$.
\begin{theorem}
\label{theorem push forward}
Suppose we have a Cartesian diagram such that
\begin{equation}
\begin{diagram}
  \node{F_2 } \arrow{e,t}{f} \arrow{s,l}{p_2} \node{F_1} \arrow{s,r}{p_1} \\
  \node{V_2} \arrow{e,b}{g} \node{V_1}
\end{diagram}
\label{eq virtual square}
\end{equation}
that
\begin{itemize}
  \item $F_i$ are Deligne-Mumford stacks.
  \item $V_i$ are Artin stacks of the same pure dimension.
  \item $g$ is a morphism of relative Deligne-Mumford type, and of pure degree $d$ for
some $d \in \Q_{\ge 0}$.
\item
$f$ is proper.
\item
$F_1 \to V_1$ has perfect relative obstruction theory $E_1$, inducing a perfect
relative obstruction theory $E_2 = f^\ast E_1$ on $F_2 \to V_2$.
\end{itemize}
Then
\begin{equation*}
  f_\ast[F_2, E_2] = d[F_1, E_1]
\end{equation*}
\end{theorem}
\begin{proof}
Let $\mf C_{F_i/V_i}$ be the relative intrinsic normal cone stack. First, we
reduce to proving that $\mf C_{F_2/V_2} \to \mf C_{F_1/V_1}$ is  of pure degree
$d$ (observe that $\mf C_{F_i/V_i}$ are of the same pure dimension). As, if
$E_1 = E^{-1}_1 \to E^0_1$, let $E_{1,1} = {E_1^{-1}}^\vee$ and let $E_{1,0} =
{E_1^0}^\vee$. Recall we have a closed embedding $\mf C_{F_1/V_1} \into
[E_{1,1}/E_{1,0}]$, where this is the stack quotient.  Let $C_1 = \mf
C_{F_1/V_1} \times _{[E_{1,1}/E_{1,0}]} E_{1,1}$. $C_1$ is an ordinary cone, so
$C_1 \to F_1$ is a scheme, and there is a closed embedding $C_1 \into E_{1,1}$.
In a similar way define $C_2 \into E_{2,1}$. The map $C_1 \to \mf C_{F_1/V_1}$
is smooth and surjective, and $C_2 = \mf C_{F_2/V_2} \times_{\mf C_{F_1/V_1} }
C_1$.  $\mf C_{F_2/V_2} \to \mf C_{F_1/V_1}$ being of pure degree $d$ is
equivalent to $C_2 \to C_1$ being of pure degree $d$. Form the diagram
\begin{equation*}
\begin{diagram}
  \node{C_2} \arrow{e} \arrow{s} \node{C_1} \arrow{s} \\
    \node{E_{2,1}} \arrow{e,t}{f'} \arrow{s,l}{p_2} \node{E_{1,1}} \arrow{s,r}{p_1} \\
  \node{F_2} \arrow{e,t}{f} \node{F_1}
\end{diagram}
\end{equation*}
By definition, $p_i^\ast [F_i, E_i] = [C_i]$ in $A_\ast E_{i,1}$. To show
$f_\ast[F_2,E_2] = d[F_1,E_1]$ is equivalent to showing $f'_\ast[C_2] = d
[C_1]$, for which it is enough to show that $C_2 \to C_1$ is of pure degree
$d$.

To do this we work locally, and reduce to the case of schemes. We pick a local
chart $U_1 \to V_1$, where $U_1$ is an irreducible scheme, and $U_1 \to V_1$ is
smooth.  Pick an \'{e}tale map of degree $n$, $U_2 \to U_1 \times_{V_1} V_2$,
where $U_2$ is a scheme, and an \'{e}tale map from a scheme $X_1 \to F_1$. By
possibly passing to smaller charts, we can pick a factorization of the map $X_1
\to U_1$ into $X_1 \into M_1 \to U_1$, where $M_1$ is a scheme, $X_1 \into M_1$
is a closed embedding, and $M_1 \to U_1$ is smooth. Without loss of generality,
we can assume that $U_2$, $X_1$ and $M_1$ are irreducible. We have a diagram
\begin{equation*}
\begin{diagram}
  \node{X_2} \arrow{e,t}{f_X} \arrow{s,l}{i_2} \node{X_1} \arrow{s,r}{i_1} \\
  \node{M_2} \arrow{e,t}{f_M} \arrow{s,l}{p_2} \node{M_1} \arrow{s,r}{p_1} \\
  \node{U_2} \arrow{e,t}{f_U} \node{U_1}
\end{diagram}
\end{equation*}
$U_2 \to U_1$ and $M_2 \to M_1$ are of degree $dn$. It is sufficient to show
that $\mf C_{X_2/U_2} \to \mf C_{X_1/U_1}$ is of pure degree $dn$. $\mf
C_{X_i/U_i} = [C_{X_i/M_i} / i_i^\ast T_{M_i/U_i}]$, where $C_{X_i/M_i}$ is the
usual normal cone. Since $f_M^\ast T_{M_1} = T_{M_2}$, it is sufficient to show
that the map $C_{X_2/M_2} \to C_{X_1/M_1}$ is of pure degree $dn$, and we have
reduced to the case of closed embeddings of schemes.

Now, we prove it in this case using the flat deformation to the normal cone, as
in \cite{fulton}, and the fact that the degree is constant in a flat family of
maps. Let $Z_i$ be the blowup of $M_i \times \mbb P^1$ along $X_i \times \{
\infty \}$. Let $M_i'$ be the blowup of $M_i$ along $X_i$. There is a
commutative diagram
\begin{equation*}
\begin{diagram}
  \node{Z_2} \arrow{e,t}{g} \arrow{se,b}{f_2} \node{Z_1} \arrow{s,r}{f_1} \\
  \node[2]{\mbb P^1}
\end{diagram}
\end{equation*}
The maps  $f_i : Z_i \to \mbb P^1$ are flat, and $f_i^{-1}(\infty) = \mbb P(C_i
\oplus 1) + M_i'$, as Cartier divisors, with multiplicity.  Clearly $Z_2 \to
Z_1$ is of pure degree $dn$, as is $M_2' \to M_1'$.  It follows that  the map
$\mbb P(C_2 \oplus 1) \to \mbb P(C_1 \oplus 1)$ is of pure degree $dn$. There
are open embeddings $C_i \into \mbb P(C_i \oplus 1)$, which implies $C_2 \to
C_1$ is of pure degree $dn$ as desired.
\end{proof}

\section{Finite degree theorem}
\label{sec_push_forward}
Let $v \in \rvertex$.  We want to construct $\eta \in \cvertex$, together with
a set of tails $A \subset T(s(\eta))$, such that:
\begin{itemize}
  \item $s(\eta)$ has just one  vertex.
  \item $v$ is obtained from $r(s(\eta))$ by removing the tails $A$.
  \item $g(t(\eta)) = 0$.
  \item $\dim \mod_\eta = \dim \mod_v$ and the map $\mod_\eta \to \mod_v$ is
  of degree $d \in \Q_{> 0}$.
\end{itemize}
The idea is quite simple: if $C \in \mod_v$ is a generic genus $g(v)$ curve
with some marked points $P_i$, and $D = \sum \lambda_i P_i$ is a positive
divisor of degree $g+1$, there is a unique (up to isomorphism) map $f :C \to
\mbb P^1$ with $f^{-1}(\infty) = D$. If everything is generic this map is
simply ramified. There is a unique \'{e}tale map of twisted curves $\mc C \to
\tilde{\mc C}$, such that the map $\tilde{\mc C} \to BS_{g+1}$ is
representable, which yields the ramified map $f: C \to \mbb P^1$ after taking
coarse moduli spaces. We let $\eta$ be the combinatorial data which labels this
\'{e}tale map $\mc C \to \tilde{\mc C}$, together with the marked points. There
is a map $\mod_\eta \to \mod_v$, which is generically finite.

Let us construct $\eta$ more formally.  We  assume that $g(v) > 0$  and $T(v)
\neq \emptyset$. Pick a partition $T(v) = I \coprod J$ into to subsets, and a
multiplicity function $d: I \to \Z_{> 0}$, with $\sum_{i \in I} d(i) = g+1$. We
define the covering $\eta = s(\eta) \cover t(\eta)$, by
\begin{itemize}
  \item $s(\eta)$, $t(\eta)$ have just one vertex.
  \item $g(t(\eta)) = 0$ and $g(s(\eta)) = g$.
  \item $a(s(\eta)) = a(t(\eta)) = a(v) \in A$.
  \item For $k \in \Z_{> 0}$ let $[k]$ be the finite set $\{ 1, \ldots,
k\}$. Then define,
\begin{equation*}
  T(t(\eta)) = J \coprod\{ \infty\} \coprod [k]
\end{equation*}
where
\begin{equation*}
  k \defeq \# I + 3g - 1
\end{equation*}
\item
The degree of $\eta$ is $g+1$.
\item
 The tails of $s(\eta)$ are
\begin{equation*}
  T(s(\eta)) =  (J\times [ g+1]) \coprod I \coprod ([g]
\times [k ])
\end{equation*}
\item
The map $T(s(\eta)) \to T(t(\eta))$ sends $I \to \infty$, and is the natural
product map on the other factors,
\begin{align*}
  J\times [g+1] & \to J \\
  [g] \times [k] & \to [k]
\end{align*}
\item
We define the multiplicity function $m$ on $T(t(\eta))$. $m(\infty)$ is the
lowest common multiple of $d(i)$ for $i \in I$,  $m(j) = 1$ for $j \in J$ and
$m(r) = 2$ for $r \in [k]$.
\item
We define the multiplicity function on $s(\eta)$. $m(i) = m(\infty)/d(i)$ for
$i \in I$, $m(r,j) = 1$ for $(r,j) \in [g+1] \times J$, and on $[g] \times[k]$
$m$ is defined by
\begin{align*}
    m(1,s) & = 1 \\
    m(r,s) & = 2\text{ if } r>1.
\end{align*}
This implies that for $(r,s) \in [g] \times [k]$,
\begin{align*}
    d(1,s) & = 2 \\
    d(r,s) & = 1\text{ if } r>1.
\end{align*}
\end{itemize}

The Riemann-Hurwitz formula becomes
\begin{equation*}
  2g - 2 = -2(g+1) + g+1 - \# I + k
\end{equation*}
which is true by our choice of $k$.  The formulae
 for the dimensions
of $\mod_\eta$ and $\mod_v$ are
\begin{align*}
 \op{dim}  \mod_\eta & = k + \# J - 2 \\
 \op{dim} \mod_v & = 3g - 3 + \# I + \# J
\end{align*}
which are  equal.

The map $\mod_\eta \to \mod_v$, comes from forgetting the tails
\begin{equation*}
  (J \times [g]) \coprod ([g] \times[k]) \into (J \times [g+1]) \coprod ([g] \times [k]) \coprod I =
  T(s(\eta))
\end{equation*}
\begin{lemma}
The map $\mod_\eta \to \mod_v$ is of degree
\begin{equation*}
 \frac{ k! (g!)^ {\# J} ((g-1)!)^{k}}{2^k m(\infty) }
\end{equation*}
\end{lemma}
\begin{proof}
Let $C \in \mod_v$ be generic, and define a divisor $D = \sum_{i \in I} d(i) i
\subset C$.  $\deg D = g+1$ and $D > 0$. For a generic curve with generic
marked points $C$, I claim that there is a unique up to isomorphism map $f : C
\to \mbb P^1$ with $f^{-1}(\infty) = D$, and further $f$ is simply ramified.

As, let $D'$ be a divisor on $C$ with $0 \le D' < D$. Let $I' \subset I$ be the
set of points which occur with non-zero multiplicity in $D'$.  Firstly, we
would like to show that the locus of smooth curves $C$ which admit a map $f : C
\to \mbb P^1$ with $f^{-1}(\infty) = D'$ is of positive codimension in
$\mod_v$. Any such curve $C \in \mod_v$, with its marked points, is determined
up to finite ambiguity by the branch points in $\mbb P^1$ of the map $f : C \to
\mbb P^1$, up to $\C \ltimes \C^\ast$-action, and by the marked points $J
\coprod (I \setminus I') \subset C$. It follows from the Riemann-Hurwitz
formula that $f$ has strictly less than $3g - 1 + \#I'$ branch points. The $\C
\ltimes \C^\ast$-action on $\mbb P^1$ reduces the dimension of the space of
possible branch points by $2$; so we find that the moduli space of smooth
curves $C$ with marked points, which admit such a map $f: C \to \mbb P^1$, is
of dimension strictly less than $3g -3 + \#J + \#I$. That is, it is of positive
codimension in $\mod_v$.

This implies that for a generic curve $C$, for all $0 \le D' < D$, $\op{dim
}H^0 (C, \Oo_C(D')) = 1$. Riemann-Roch tells us that $\op{dim}H^0(C, \Oo_C(D))
\ge 2$. Take $D'$ to be of degree $g$. We must have $H^1(C, \Oo_C(D')) = 0$,
which implies $H^1(C, \Oo_C(D)) = 0$ and $\op{dim}H^0(C, \Oo_C(D)) = 2$.  This
last fact implies that there is precisely one map $f: C \to \mbb P^1$, up to
isomorphism, such that $f^{-1}(\infty) = D$. I claim that if $C$ is generic,
this map is simply ramified.  One can see this by observing, using the
Riemann-Hurwitz formula as before, that the locus of smooth curves $C \in
\mod_v$ which admit a map $f : C \to \mbb P^1$ with $f^{-1}(\infty) = D$ and
non-simple ramification is of positive codimension.

The degree of our map $\mod_\eta \to \mod_v$ can now be calculated from
different ways of ordering tails of $\eta$, and automorphisms of twisted curves
in $\mod_\eta$ over their coarse moduli space.
\end{proof}

\section{Stable curves in $X$}
\label{sec_base_change}

Let $X$ be a smooth projective variety.  Let $C(X)$ be the Mori cone of curves
in $X$ modulo numerical equivalence.  $C(X)$ is a semigroup with indecomposable
zero and finite decomposition.  We define our categories of graphs and vertices
using $C(X)$.

For every $\gamma \in \rgraph$, we have \cite{behrend gw} the moduli stack
$\cmod_\gamma(X)$ of stable maps to $X$ of type $\gamma$. $\cmod_\gamma(X)$ is
a separated, proper Deligne-Mumford stack of finite type. There is a map
$\cmod_\gamma(X) \to \mod_\gamma$, with a perfect relative obstruction theory
$(R\pi_\ast f^\ast TX)^\vee$, where $\pi: \mc C \to \cmod_\gamma(X)$ is the
universal curve and $f: \mc C \to X$ the universal map.  The target
$\mod_\gamma$ for this perfect relative obstruction theory is slightly
different to the version used in Behrend's construction \cite{behrend gw},
because of the labellings by elements of the semigroup $A$.  The virtual
fundamental classes, however, are the same.  This follows from the fact that
the map $\mod_{g,n,a} \to \mod_{g,n}$ is \'etale.

Suppose we have a map $\gamma' \to \gamma$ in $\rgraph$. Then we have a fibre
square,
\begin{equation*}
\begin{diagram}
  \node{\cmod_{\gamma'}(X)} \arrow{e} \arrow{s,r}{p'} \node{\cmod_\gamma(X)} \arrow{s,l}{p} \\
  \node{\mod_{\gamma'}} \arrow{e} \node{\mod_\gamma}
\end{diagram}
\end{equation*}
Further, the perfect relative obstruction theory of $p'$ is pulled back from
that of $p$.

If $\gamma'$ is obtained from $\gamma$ by adding on some tails, then  we have a
fibre square exactly as above, and again the perfect relative obstruction
theory of $p': \cmod_{\gamma'}(X) \to \mod_{\gamma'}$ is pulled back from that
of $p: \cmod_{\gamma}(X) \to \mod_\gamma$.

If $\gamma'$ is obtained by cutting an edge of $\gamma$, then $\mod_{\gamma'} =
\mod_\gamma$.  We have a fibre square,
\begin{equation*}
\begin{diagram}
  \node{\cmod_{\gamma}(X)} \arrow{e} \arrow{s} \node{\cmod_{\gamma'}} \arrow{s} \\
  \node{X \times \mod_{\gamma}} \arrow{e,t}{\bigtriangleup} \node{X\times X \times \mod_{\gamma'}}
\end{diagram}
\end{equation*}
The perfect relative obstruction theories of $\cmod_{\gamma'}(X)$ and
$\cmod_{\gamma}(X)$ over $\mod_{\gamma'} = \mod_\gamma$ are compatible with
this Cartesian diagram, in the sense of \cite{bf}, section 7.

Let $\eta \in \cov$.  Define $\cmod_\eta(X)$ by the Cartesian square
\begin{equation*}
\begin{diagram}
  \node{\cmod_\eta(X)} \arrow{e,t}{(r\circ s)_X} \arrow{s} \node{\cmod_{r(s(\eta))}} \arrow{s} \\
  \node{\mod_\eta} \arrow{e,t}{r \circ s} \node{\mod_{r(s(\eta))}}
\end{diagram}
\end{equation*}
This is the stack of diagrams $\mc C \from \mc C' \to X$, where $\mc C' \to \mc
C$ is a map from $\mod_\eta$ and, if $C'$ is the coarse moduli space of $\mc
C'$, the map $C' \to X$ is a stable map from $\cmod_{r(s(\eta))}(X)$.  Give
$\cmod_\eta(X) \to \mod_\eta$ the perfect relative obstruction theory pulled
back from that for $\cmod_{r(s(\eta))}(X) \to \mod_{r(s(\eta))}$.

If $\eta' \to \eta$ is a map in $\cov$, then we have a fibre square
\begin{equation*}
\begin{diagram}
  \node{\cmod_{\eta'}(X)} \arrow{e} \arrow{s,r}{p'} \node{\cmod_\eta(X)} \arrow{s,l}{p} \\
  \node{\mod_{\eta'}} \arrow{e} \node{\mod_\eta}
\end{diagram}
\end{equation*}
and the perfect relative obstruction theory for $p'$ is pulled back from that
for $p$.

Let $\eta \in \cov$ and let $e \in E(t(\eta))$. Let $I \subset E(s(\eta))$ be
the set of edges lying over $e$.  Let $\eta' \in \cov$ be obtained from $\eta$
by cutting the edges $e,I$.  Then, $\mod_{\eta'} = \mod_\eta$, and we have a
Cartesian diagram
\begin{equation*}
\begin{diagram}
  \node{\cmod_{\eta}} \arrow{e} \arrow{s} \node{\cmod_{\eta'}} \arrow{s} \\
  \node{\mod_{\eta} \times X^I} \arrow{e,t}{\diag} \node{\mod_{\eta'} \times (X^2)^I}
\end{diagram}
\end{equation*}
which is compatible with perfect relative obstruction theories over $\mod_\eta
= \mod_{\eta'}$.

Observe that since $\mod_\eta \to \mod_{t(\eta)}$ is \'{e}tale, $\cmod_\eta(X)$
has a perfect relative obstruction theory over $\mod_{t(\eta)}$ also.
\subsection{Stable maps to symmetric products}
We have described stacks of stable maps to a smooth projective variety $X$. We
also have stacks of stable maps to a smooth DM stack $V$, as defined by
Abramovich and Vistoli in \cite{av}. We are interested in these  when $V = S^d
X$. We need something to play the role of the Mori cone, that is to hold
homology classes of curves. We simply use again $C(X)$, the Mori cone of $X$.
For $\gamma \in \graph(C(X)) = \graph$, let $\cmod_\gamma (S^d X)$ be the stack
of stable representable maps from curves in $\mod_\gamma$ to $S^d X$, in a way
compatible with the $C(X)$-markings on the curve.
\begin{lemma}
For a covering $\gamma' \cover \gamma$, let $\op{Aut}(\gamma' \mid \gamma)$ be
the  group of automorphisms $\gamma'$, commuting with the covering $\gamma'
\cover \gamma$. Then,
\begin{equation*}
  \cmod_\gamma(S^d X) = \coprod_{\gamma' \cover \gamma} \cmod_{\gamma' \cover
  \gamma}(X)/ \op{Aut}(\gamma' \mid \gamma)
\end{equation*}
where the union is over all $\gamma' \cover \gamma$ of degree $d$.  Further,
this identification is compatible with perfect relative obstruction theories
over $\mod_\gamma$.
\end{lemma}
\begin{proof}
The isomorphism at the level of stacks follows from section \ref{sec map sym
prod}. We need to prove compatibility of perfect relative obstruction theories.
Suppose we have a representable stable map $f: \mc C \to S^d X$, corresponding
to a diagram $\mc C \overset{p'}{\from} \mc C' \overset{f'}{\to} X$, and
equivalently to a principal $S_d$ bundle $p : P \to \mc C$, where $P$ is an
algebraic space, and an $S_d$-equivariant map $g : P \to X^d$. The perfect
relative obstruction theory for stable maps to $S^d X$, is given by $H^\ast
(\mc C, f^\ast TS^d X)$.  But,
\begin{equation*}
  f^\ast TS^d X = p_\ast^{S_d} g^\ast T X^d = p'_\ast f'^\ast TX
\end{equation*}
Observe $p'_\ast$ and $p_\ast^{S_d}$ are exact. Let $g':C' \to X$ be the map
from the coarse moduli space of $\mc C'$ to $X$, and let $m: \mc C' \to C'$ be
the canonical map. Observe $m_\ast$ is exact, and  $f' = g' \circ m: \mc C' \to
X$. Now,
\begin{align*}
  H^\ast (\mc C, f^\ast TS^d X) &= H^\ast (\mc C', f'^\ast TX)\\
  &= H^\ast( \mc C', m^\ast g'^\ast TX) \\
  &= H^\ast( C', g'^\ast TX)
\end{align*}
as desired.
\end{proof}
There are evaluation maps $\cmod_{\gamma' \cover \gamma} (X) \to
X^{T(\gamma')}$. These come from the evaluation maps
\begin{equation*}
  \cmod_{\gamma}(S^d X)  \to \text{ twisted sectors of } S^d X
\end{equation*}
which are used to define quantum cohomology of $S^d X$. The stack of twisted
sectors $\tilde V$ of a DM stack $V$ is the stack of cyclic gerbes in $V$
\cite{agv}, \cite{chen ruan},\cite{toen riemann roch},
\begin{equation*}
  \tilde V = \coprod_{k \ge 1} \op{HomRep}( B\mu_k, V)
\end{equation*}
We can identify $\widetilde{S^d X}$  with a disjoint union $\coprod_{\sigma \in
{S_{d}}_\ast} ((X^d)^\sigma)/C(\sigma)$, where the disjoint union is over
conjugacy classes in $S_d$, $\sigma$ is a representative of each conjugacy
class, $(X^d)^\sigma$ is the $\sigma$-fixed points and $C(\sigma)$ is the
centralizer of $\sigma$. There is a commutative diagram of evaluation maps
\begin{equation*}
\begin{diagram}
  \node{\cmod_{\gamma' \cover \gamma}(X)} \arrow{e,t}{ev} \arrow{s} \node{X^{T(\gamma')}} \arrow{s} \\
  \node{\cmod_\gamma(S^d X)} \arrow{e,t}{ev} \node{(\widetilde{S^d X})^{T(\gamma)}}
\end{diagram}
\end{equation*}

Further, the tautological $\psi$-classes on $\cmod_\gamma(S^d X)$ are pulled
back to $\psi$ classes on $\cmod_{\gamma' \cover \gamma}(X)$. Thus one can
identify integrals of the form
\begin{equation*}
  \int_{[\cmod_{\gamma' \cover \gamma}]_{virt}} \prod_{t \in T(\gamma)}
  \psi_t^{k_t} \prod_{t' \in T(\gamma')} ev_{t'}^\ast h_{t'}
\end{equation*}
where $h_{t'} \in H^\ast (X)$, with Gromov-Witten invariants of $S^d X$.

\section{From genus $g$ invariants of $X$ to genus $0$ invariants of $S^d X$}
\label{sec main theorem}
Let $v \in \rvertex$, so that $v$ labels a stack of stable maps to $X$. Assume
$g(v) > 0$ and $\#T(v) > 0$. We will use the notation of section
\ref{sec_push_forward}. There we constructed $\eta \in \cvertex$, such that $v$
is obtained by removing some tails of $r(s(\eta))$, and $g(t(\eta)) = 0$.  The
associated map
\begin{equation*}
  \mod_\eta \to \mod_v
\end{equation*}
was shown to be of degree
\begin{equation*}
 n \defeq \frac{k! (g!)^ {\# J} ((g-1)!)^{k}}{2^k m(\infty) }
\end{equation*}
Form the fibre square
\begin{equation*}
\begin{diagram}
  \node{\cmod_\eta(X)} \arrow{e,t}{q} \arrow{s} \node{\cmod_v(X)} \arrow{s} \\
  \node{\mod_\eta} \arrow{e} \node{\mod_v}
\end{diagram}
\end{equation*}
\begin{lemma}
The map $q: \cmod_\eta(X) \to \cmod_v(X)$, is of degree $n$, in the virtual
sense,
\begin{equation*}
  q_\ast [\cmod_\eta(X)]_{virt} = n [\cmod_v(X)]_{virt}
\end{equation*}
\end{lemma}
\begin{proof}
We apply theorem \ref{theorem push forward}. Observe that $\mod_\eta \to
\mod_v$ is  relatively of Deligne-Mumford type and generically finite of degree
$n$,  $\mod_\eta$ and $\mod_v$ are algebraic stacks, and that $\cmod_\eta(X)
\to \cmod_v(X)$ is proper.
\end{proof}
We have seen in section \ref{sec_pull_back} how to express the pulled-back
tautological classes $p^\ast \psi_t$, and their products, for $ t\in T(v)$, in
terms of tautological classes pushed forward under contractions $\rho \to \eta$
in $\cov$. Let us combine this result with the previous one to calculate
Gromov-Witten invariants of $X$ in terms of integrals over $\cmod_\rho(X)$,
which are genus $0$ invariants of $S^d X$.

We have a commutative diagram
\begin{equation}
\begin{diagram}
  \node{X^{T(\eta)}} \arrow{e,t}{\pi}  \node{X^{T(v)}} \\
  \node{\cmod_\eta(X)} \arrow{e,t}{q}\arrow{n,l}{ev_\eta} \arrow{s,r}{c_\eta} \node{\cmod_v(X)} \arrow{n,r}{ev_v} \arrow{s,l}{c_v} \\
  \node{\mod_\eta} \arrow{e,b}{p} \node{\mod_v}
\end{diagram}
\label{eq diagram main theorem}
\end{equation}
The integrals we want to calculate are
\begin{equation*}
  \int_{[\cmod_v(X)]_{virt}} c_v^\ast e^{\sum_{t \in T(v)} z_t \psi_t} ev_v^\ast \alpha
\end{equation*}
where $\alpha \in H^\ast X^{T(v)}$.

Let us recall some of the notation of section \ref{pull back psi classes}.  For
each map $f: \rho \to \eta$ in $\cov$, we defined
\begin{equation*}
  \norm f = \frac{\#\op{Aut}(\rho \to \eta \mid \eta)}{ \prod_{e \in E(t(\rho))} m(e)^2}
\end{equation*}
where $\op{Aut}(\rho \to \eta \mid \eta)$ is the group of automorphisms of
$\rho$ commuting with the contraction $\rho \to \eta$, or equivalently (in this
special case) fixing all tails.

Let $I \subset T(s(\eta))$ be the set of tails we forget to obtain $v$. For
each $f: \rho \to \eta$, each edge $e \in E(s(\rho))$, and each tail $t \in v$,
we defined $S(e,t,I) \in \{0,1\}$ as follows.  Let $s(\rho)_e$ be obtained from
contracting all edges other than $e$. If the vertex of $s(\rho)_e$ containing
$t$ becomes unstable after forgetting the tails $I$, we set $S(e,t,I) = 1$,
otherwise $S(e,t,I) = 0$.

For each tail $t \in T(v)$, define a formal variable $z_t$, and for each edge
$e \in E(t(\rho))$, define a variable $w_e$, with the relations
\begin{equation*}
  w_e = \sum_{e' \in p^{-1}(e) \subset E(s(\rho))} m(e') S(e',t,I) z_t
\end{equation*}
Using this notation, we have
\begin{theorem}[Main theorem]
\begin{multline}
\label{main formula}
\int_{[\cmod_v(X)]_{virt}} c_v^\ast e^{\sum_{t \in T(v)} z_t \psi_t} ev_v^\ast
\alpha = \\
\sum_{f: \rho \to \eta} \tfrac{1}{n \norm f}\int_{[\cmod_\rho(X)]_{virt}}
c_\rho^\ast\left(  e^{\sum_{t
 \in T(t(\rho))} w_t \psi_t} \prod_{e \in E(t(\rho))} \frac{ 1 - e^{w_e \psi_e}  }{\psi_e}\right) ev_\rho^\ast
 \pi^\ast \alpha
\end{multline}
The left hand side is the general form for descendent genus $g$ invariants of
$X$. The right hand side is an expression in the genus $0$ invariants of the
symmetric product stack  $S^{g+1}X$.
\end{theorem}
\begin{proof}
Firstly, the projection formula shows that
\begin{equation*}
  \int_{[\cmod_v(X)]_{virt}}  c_v^\ast e^{\sum_{t \in T(v)} z_t \psi_t} ev_v^\ast \alpha
 = \tfrac{ 1}{n} \int_{[\cmod_\eta(X)]_{virt}} q^\ast c_v^\ast e^{\sum_{t \in T(v)} z_t
 \psi_t}q^\ast ev_v ^\ast \alpha
\end{equation*}
The formula \ref{eq psi pullback} shows that
\begin{equation*}
  q^\ast c_v^\ast e^{\sum_{t \in T(v)} z_t \psi_t} = c_\eta^\ast \left( \sum_{f: \rho \to \eta} e^{\sum_{t
 \in T(t(\rho))} w_t \psi_t} \prod_{e \in E(t(\rho))} \frac{ 1 - e^{w_e \psi_e} }{\psi_e}
 c_\eta^\ast [\mod_f]
 \right)
\end{equation*}
 Now, for each term
in the sum, $\mod_\rho \to \mod_f$ is \'{e}tale of degree $\norm f$.  The
standard compatibility of virtual fundamental classes shows that
\begin{equation*}
  c_\eta^\ast [\mod_f] = \tfrac{1}{\norm{f}} f_\ast [\cmod_\rho (X)]
\end{equation*}
where $f_\ast: \cmod_\rho (X) \to \cmod_\eta(X)$ is the canonical map. This
implies the result.
\end{proof}

\section{Examples}
\label{section calculations}
\subsection{Generalities}
I will calculate some examples in the case where $X$ is a point.  We will work
with the semigroup $A = 0$. For each $\eta \in \cvertex$, we have the
associated moduli stack $\cmod_\eta$, with a map
\begin{equation*}
 p: \cmod_\eta \to \cmod_{g(t(\eta)), T(t(\eta))}
\end{equation*}
The stack $\cmod_{g(t(\eta)), T(t(\eta))}$ is the usual Deligne-Mumford stack
of stable curves. It follows from the results of \cite{acv} that $\cmod_\eta$
is closely related to the normalization of a stack of admissible covers.

In this section, we will always assume $g(t(\eta)) = 0$, and that $\#T(t(\eta))
= n$. Further, we pick an ordering on the set $T(t(\eta))$, and so an
isomorphism
\begin{equation*}
  T(t(\eta)) \iso [n] = \{1, \ldots, n\}
\end{equation*}
The map $p$ is now a map $\cmod_\eta \to \cmod_{0,n}$; we want to calculate its
degree. For each vertex $v \in V(s(\eta))$, and tail $t' \in T(s(\eta))$,
recall we have numbers $d(v), d(t') \in \Z_{\ge 1}$. Let $i \in [n] \iso
T(t(\eta))$. Then, the set
\begin{equation*}
  \{ d(t') \mid t' \in p^{-1}(i) \subset T(s(\eta)), \, t' \text{ is attached
  to } v \}
\end{equation*}
defines  a partition of $d(v)$, and so a conjugacy class $C_{i,v} \subset
S_{d(v)}$.  Let
\begin{equation*}
  \chi(v) = \{ \sigma_i  \in S_{d(v)} \text{ for  } i = 1 \ldots n  \mid
  \sigma_i \in C_{i,v}, \, \prod_{i = 1}^n \sigma_i = 1, \, \#( [d(v)] / \ip{\sigma_1, \ldots ,\sigma_n}) = 1 \}
\end{equation*}
The last condition means that the group $\ip{\sigma_1, \ldots, \sigma_n}$ acts
transitively on the set $[d(v)] = \{ 1, \ldots, d(v) \}$.

Let $\op{Aut}(\eta \mid t(\eta), V(s(\eta)))$ be the group of automorphisms of
$\eta$ acting trivially on $t(\eta)$ and the set of vertices $V(s(\eta))$.
\begin{proposition}
The map $\cmod_\eta \to \cmod_{0,n}$ is of degree
\begin{equation*}
  \frac{\#\op{Aut}(\eta \mid t(\eta), V(s(\eta)))} {\prod_{i = 1}^n m(i)}
  \prod_{v \in V(s(\eta))}
  \frac{ \chi (v) }{d(v)!}
\end{equation*}
\end{proposition}
Recall that for each $i \in [n]$,
\begin{equation*}
  p^\ast \psi_i = m(i) \psi_i
\end{equation*}
Let $k_i \in \Z_{\ge 0}$, for $i = 1 \ldots n$. We have
\begin{equation}
\label{formula degree}
  \int_{\cmod_\eta} \prod_{i = 1}^n \psi_i^{k_i} =
  \frac{\#\op{Aut}(\eta \mid t(\eta), V(s(\eta)))} {\prod_{i = 1}^n m(i)^{k_i + 1}} \prod_{v \in V(s(\eta))} \frac{
  \chi(v)}{d(v)!} \int_{\cmod_{0,n}} \prod_{i = 1}^n \psi_i^{k_i}
\end{equation}
which allows us, at least in principle, to calculate the integrals on the left
hand side.

\subsection{Calculations}
The first example we will compute is
\begin{example}
\begin{equation*}
  \int_{\cmod_{1,1}} \psi_1 = 1/24
\end{equation*}
\end{example}
We define $\eta$ by
\begin{itemize}
  \item $s(\eta)$ has just one vertex; $g(s(\eta)) = 1$ and
  $g(t(\eta)) = 0$.
  \item The tails are
\begin{align*}
  T(s(\eta)) &= \{ X_1, X_2, X_3, X_4 \} & d(X_i) &= 2 \quad m(X_i) = 1 \\
  T(t(\eta)) &= \{x_1,x_2,x_3,x_4\} & m(x_i) &= 2
\end{align*}
  \item The map $T(s(\eta)) \to T(t(\eta))$ sends $X_i \mapsto x_i$.
\end{itemize}
Forgetting the marked points $X_2,X_3,X_4$ gives us a map
\begin{equation*}
  \pi : \cmod_\eta \to \cmod_{1,1}
\end{equation*}
Now we apply the main theorem.  Note that the only graph $f: \rho \to \eta$
that arises with non-zero coefficient on the right hand side of (\ref{main
formula}) is $\rho = \eta$. Therefore
\begin{equation*}
  \int_{\cmod_{1,1}} \psi_1 = 3^{-1} \cdot 2^3 \int_{\cmod_\eta} \psi_1
\end{equation*}
Next, we apply formula (\ref{formula degree}). In this case, $\chi(v) = 1$ for
the unique vertex $v \in V(s(\eta))$, and $\# \op{Aut}( \eta \mid t(\eta),
V(s(\eta))) = 1$, so we find that
\begin{equation*}
  \int_{\cmod_\eta} \psi_1 = 2^{-6} \int_{\cmod_{0,4}} \psi_1 = 2^{-6}
\end{equation*}
Combining these formulae yields
\begin{equation*}
  \int_{\cmod_{1,1}} \psi_1 = 1/24
\end{equation*}
Next we calculate
\begin{example}
\begin{equation*}
  \int_{\cmod_{2,1}} \psi^4 = 1/1152 = 2^{-7} 3^{-2}
\end{equation*}
\end{example}
One can see this is correct by applying the Kontsevich-Witten theorem. We
define $\eta \in \cvertex$ in this case by
\begin{itemize}
  \item $s(\eta)$ has just one vertex, with $g(s(\eta)) =
  2$ and $g(t(\eta)) = 0$.
  \item We define the sets of tails by
\begin{align*}
  T(s(\eta)) &= \{ X_1, X_2, Y_2, X_3, Y_3, \ldots,  X_7, Y_7 \} \\
  T(t(\eta)) &= \{x_1,x_2,x_3, \ldots, x_7\}
\end{align*}
with the map $T(s(\eta)) \to T(t(\eta))$, sending $X_i \mapsto x_i$ and $Y_j
\mapsto x_j$.
\item The multiplicities of these tails are defined by
\begin{align*}
  d(X_1) &= 3 & d(X_i) &= 2 \text{ if } i \ge 2 \\
  m(X_i) &= 1 \text{ for all } i\\
  d(Y_i) &= 1 & m(Y_i) &= 2 \\
  m(x_1) &= 3 & m(x_i) &= 2 \text{ if } i \ge 2
\end{align*}
\end{itemize}
Forgetting the marked points $X_i,Y_i$ for all $i \ge 2$, gives a map
\begin{equation*}
  \pi : \cmod_\eta \to \cmod_{2,1}
\end{equation*}
Now we apply the main formula (\ref{main formula}).  In this case, we find that
there are non-trivial graphs $f : \rho \to \eta$ appearing on the right hand
side.  For each $i = 2 \ldots 7$ define a graph $\rho_i \in \cov$ with a map
$f_i : \rho_i \to \eta$, as follows.
\begin{itemize}
  \item $s(\rho_i)$ has 3 vertices, $V_1,V_2,W_2$, and $t(\rho_i)$ has two
  vertices $v_1, v_2$.  The map $V(s(\rho_i)) \to V(t(\rho_i))$ sends $V_i \to
  v_i$ and $W_2 \to v_2$.
  \item The genera of the vertices are given by
\begin{equation*}
  g(v_i) = 0 \quad g(V_1) = g(W_2) = 0 \quad g(V_2) = 2
\end{equation*}
  \item Their is an edge $e$ joining $v_1$ and $v_2$, an edge $E$ joining
  $V_1$ and $V_2$, and an edge $F$ joining $V_1$ and $W_2$.  The multiplicities
  of these edges is given by
\begin{align*}
  m(e) &= 2 & m(E) &= 1 & m(F) &= 2\\
  &&d(E) &= 2 &d(F) &= 1
\end{align*}
\item
The sets of tails of the vertices are given by
\begin{equation*}
\begin{array}{ccc }
  T(V_1) = \{X_1,X_i,Y_i\} & T(V_2) = \{X_j \mid j \neq 1,i\} & T(W_2) = \{Y_j \mid j \neq i\} \\
  T(v_1) = \{x_1,x_i\} & T(v_2) = \{x_i \mid i \neq 1,i\} & \\
\end{array}
\end{equation*}
\end{itemize}
Observe that after contracting the edge $F$, $X_1$ is on a genus $0$ vertex
that becomes unstable after removing the other marked points.  This implies
$f_i: \rho_i \to \eta$ occurs with non-zero coefficient in the expansion in the
main formula (\ref{main formula}). In fact, $\rho_i$ are the only non-trivial
graphs which occur. We have $\frac{1}{\norm{f_i}} = 4$. Note also that on
$\cmod_{\rho_i}$, $\psi_1 = 0$ because the marked point $x_1$ is on a genus $0$
curve with two marked points and one edge.  Applying the main formula, we see
\begin{equation*}
  \int_{\cmod_{2,1}} \psi^4 = \frac{2^6 \cdot 3}{6!}
  \left( \int_{\cmod_\eta}  \psi_1^4 -  4 \sum_{i = 2}^7 \int_{\cmod_{\rho_i}} \psi_e^3 \right)
\end{equation*}

Next, we apply formula (\ref{formula degree}) to calculate $\int_{\cmod_\eta}
\psi_1^4$.  One can calculate easily that
\begin{equation*}
  \# \{ \sigma_1, \sigma_2, \ldots, \sigma_7 \in S_3 \mid \sigma_1
  \text{ is a 3-cycle}, \sigma_i \text{ are transpositions for } i \ge 2, \prod \sigma_j = 1
  \}= 3^5 \cdot 2
\end{equation*}
Further, $\#\op{Aut}(\eta \mid t(\eta), V(s(\eta))) =
1$, so that
\begin{equation*}
  \int_{\cmod_\eta}\psi_1^4 =  2^{-6} \cdot 3^{-1} \int_{\cmod_{0,7}} \psi_1^4
  = 2^{-6} \cdot 3^{-1}
\end{equation*}
Now we calculate $\int_{\cmod_{\rho_i}} \psi_e^3$.  As in subsection
\ref{subsection coverings}, $\cmod_{\rho_i}$ splits as a product of
contributions from the vertices of $t(\rho_i)$. Let $p: V(s(\rho_i)) \to
V(t(\rho_i))$ be the natural map. For each $v \in V(t(\rho_i))$, let $p^{-1}(v)
\cover v \in \cvertex$ be the natural covering, whose tails consist of tails
and germs of edges in $\rho_i$ at the vertices $v,p^{-1}(v)$. We have
\begin{equation*}
  \cmod_{\rho_i} = \prod_{v \in V(t(\rho_i))} \cmod_{p^{-1}(v) \cover v}
\end{equation*}
This implies that integrals split in a similar way.  There  are two vertices
$v_1, v_2$ on $t(\rho_i)$. We have $p^{-1}(v_1) = V_1$ and $p^{-1}(v_2) =
\{V_2,W_2\}$.  Further,
\begin{equation*}
  \op{dim} \cmod_{p^{-1}(v_1) \cover v_1} = 0
\end{equation*}
Denote by $e$ the germ of the edge at $v_2$, which we consider as a tail.  We
have
\begin{equation*}
  \int_{\cmod_{\rho_i}}\psi_e^3 = \left( \int_{\cmod_{p^{-1}(v_1) \cover v_1}} 1 \right)\times
  \left( \int_{\cmod_{p^{-1} (v_2) \cover v_2}} \psi_e^3 \right)
\end{equation*}
An easy application of formula (\ref{formula degree}) now shows that
\begin{align*}
  \int_{\cmod_{p^{-1}(v_1) \cover v_1}} 1 &= 3^{-1}\cdot 2^{-2} \\
  \int_{\cmod_{p^{-1} (v_2) \cover v_2}} \psi_e^3 & = 2^{-10}
\end{align*}
Finally, we see that
\begin{align*}
  \int_{\cmod_{2,1}} \psi^4 &= \frac{2^6 \cdot 3}{6!}\left(2^{-6} \cdot 3^{-1} - 24 \cdot 2^{-12} \cdot 3^{-1}
  \right)\\
  &= 2^{-7} \cdot 3^{-2}
\end{align*}
as desired.

\end{document}